\documentclass[oneside,english]{amsart}

\usepackage[latin1]{inputenc}
\usepackage{amssymb}

\usepackage{amssymb}

\makeatletter

\author{\'{E}lise Janvresse, Tom Meyerovitch,\\ Emmanuel Roy and Thierry de la Rue}
\title{Poisson suspensions and Entropy for infinite transformations}
\address{\'Elise Janvresse, Thierry de la Rue:
Laboratoire de Math\'ematiques Rapha\"el Salem,
Universit\'e de Rouen, CNRS --
Avenue de l'Universit\'e --
F76801 Saint \'Etienne du Rouvray, France.}
\email{Elise.Janvresse@univ-rouen.fr\\Thierry.de-la-Rue@univ-rouen.fr}
\address{Tom Meyerovitch:
School of Mathematical Sciences, Tel-Aviv University, Ramat-Aviv,
Tel-Aviv 69978, Israel}
\email{tomm@post.tau.ac.il}
\address{Emmanuel Roy: Laboratoire Analyse, G\'eom\'etrie et Applications, Universit\'e Paris 13 Institut Galil\'ee --
99 avenue Jean-Baptiste Cl\'ement --
F93430 Villetaneuse, France.}
\email{roy@math.univ-paris13.fr}

\newtheorem{theorem}{Theorem}[section]
\newtheorem{lemma}[theorem]{Lemma}
\newtheorem{proposition}[theorem]{Proposition}
\newtheorem{corollary}[theorem]{Corollary}

\usepackage{color}
\usepackage{epsfig}

\newcommand{\hide}[1]{}
\usepackage{babel}
\makeatother

\newcommand{\NN}{\mathbb{N}}

\newcommand{\B}{\mathcal{B}}

\newcommand{\hkr}{h_{\mbox{\scriptsize Kr}}} % notation for Krengel entropy
\newcommand{\hpa}{h_{\mbox{\scriptsize Pa}}} % notation for Parry entropy
\newcommand{\hrel}{h} % notation for relative entropy

\begin{document}

\subjclass[2000]{37A05, 37A35, 37A40, 28D20}
\begin{abstract}
The Poisson entropy of an infinite-measure-preserving transformation is defined in \cite{roy_thesis} as the Kolmogorov entropy of its Poisson suspension.
In this article, we relate Poisson entropy with other definitions of entropy for infinite transformations:
For quasi-finite transformations we prove that Poisson entropy coincides with Krengel's and Parry's entropy (Theorem~\ref{thm:poisson_eq_krengel_quasi_finite}).
In particular, this implies that for null-recurrent Markov chains, the usual formula for the entropy $-\sum q_i p_{i,j}\log p_{i,j}$ holds for any definitions of entropy. Poisson entropy dominates Parry's entropy  in any conservative transformation (Theorem~\ref{thm:poisson_eq_parry}). We also prove that relative entropy  (in the sense of \cite{danilenko_rudolph07}) coincides with the relative Poisson entropy (Proposition~\ref{prop:relative_entropy}).
Thus, for any factor of a conservative transformation,  difference of the Krengel's entropies equals difference of the Poisson entropies.
In case there  exists a factor with zero Poisson entropy, we prove the existence of a maximum (Pinsker) factor with zero Poisson entropy. Together with the preceding results, this answers affirmatively the question raised in \cite{aaro_park} about existence of a Pinsker factor in the sense of Krengel for  quasi-finite transformations.
\end{abstract}

\maketitle
\section{Introduction}

The basic question considered in this paper is the following:
Is there a ``natural'' entropy theory for infinite-measure-preserving dynamical systems?
Both Krengel \cite{krengel67} and Parry \cite{parry_generators69} defined notions of entropy for measure-preserving transformations,
which elegantly generalize Kolmogorov's entropy of a probability-preserving transformation. It is still an open question whether, for any conservative measure-preserving transformation, Parry's definition of entropy coincides with Krengel's. However, there are known sufficient conditions on a system for these two numbers to coincide, and Krengel's entropy dominates Parry's entropy in general.

In the present paper, we relate Parry's and Krengel's definitions of entropy with \emph{Poisson entropy}, which is the Kolmogorov entropy of the Poisson suspension, an approach previously taken in \cite{roy_thesis}. The main question here is whether Poisson entropy is equal either to Parry's entropy or to Krengel's entropy (or both) for any conservative measure-preserving transformation.
We are yet unable to completely solve this question in full generality, but we give many intermediate results. In particular, we show that Poisson entropy dominates Parry's entropy in general (Theorem~\ref{thm:poisson_eq_parry}), and that equality of all three definitions of entropy holds for many classes of transformations.

As a consequence, we obtain an intuitive expression for the entropy of Poisson-suspensions of null-recurrent Markov chains. We thus correct a mistake in \cite{grabinsky84}, where it was claimed that the entropy is always infinite for Poisson suspensions of null-recurrent Markov chains.

We prove that Poisson entropy is equal to Parry's entropy and Krengel's entropy in the following cases: quasi-finite transformations (Theorem~\ref{thm:poisson_eq_krengel_quasi_finite}) and rank-one transformations (Proposition~\ref{prop:cutting_and_stacking_zero_entropy}).
We also prove that relative entropy  (in the sense of \cite{danilenko_rudolph07}) coincides with relative Poisson entropy (Proposition~\ref{prop:relative_entropy}).
We prove that Poisson entropy is a linear functional, just as Krengel's entropy and Parry's entropy.

In section~\ref{sec:zero_entropy} we give a spectral criterion for zero Poisson entropy, which was previously shown to imply zero Parry entropy.
In Section~\ref{sec:perfect_partitions}, among other results, we prove the following dichotomy for ergodic quasi-finite infinite measure-preserving transformations: either it is remotely infinite or there exists a maximum (Pinsker) factor with zero Poisson, Krengel and Parry entropy. The proof relies on the existence of perfect Poissonian $\sigma$-algebra and illustrates the interest of using Poisson suspensions to derive results in infinite-measure-preserving ergodic theory via the far more developed finite-measure-preserving case. We also state and prove a strong disjointness result in terms of Poisson suspensions.

\textbf{Acknowledgement:} T.M. would like to thank his Ph.D. advisor, Professor Jon Aaronson, for his guidance throughout this work. The authors thank A. Danilenko and D. Rudolph for making their paper \cite{danilenko_rudolph07} available prior to its publication.

\section{Poisson suspensions and Poisson entropy}
The \emph{Poisson suspension} $(X^*,\mathcal{B}^*,\mu^*,T_*)$ of a
standard, $\sigma$-finite invertible measure-preser\-ving transformation
$(X,\mathcal{B},\mu,T)$ is a canonical method of associating a
probability-preserving transformation to a $\sigma$-finite-measure-preserving transformation.
Informally, it is a system of
non-interacting ``identical'' particles in $X$, each of which
propagates according to the transformation $T$, and such that the
expected number of particles in a set $A \in \mathcal{B}$ is
determined by $\mu(A)$.
Poisson suspensions have been studied in mathematical physics as well as in ergodic-theory
and probabilistic contexts
\cite{goldstein_lebowitz74,grabinsky84,kalikow81,MR0297288,MR0289094}
and recently in \cite{roy_thesis} and \cite{zweimuller_2007}.

There are various ways to describe a Poisson suspension. Here is one: Let $X^*$ denote the
space of measures on $X$, and let $\mathcal{B}^*$
denote the $\sigma$-algebra generated by the collection of sets
\begin{equation}
\label{eq:B_star_sigma_algebra}
\Bigl\{ \{\gamma \in X^*:~ \gamma(B) \in [a,b]\}: B \in \mathcal{B},\, 0\le a \le b \le \infty ~\Bigr\}.
\end{equation}

The probability measure $\mu^*$ on $(X^*,\mathcal{B}^*)$ is uniquely defined by requiring that measures of disjoint sets be independent and that the measure of each set $A\in\mathcal{B}$ be Poisson distributed with parameter $\mu(A)$:
$$\mu^*\left(\gamma(A)=k\right)=e^{-\mu(A)}\frac{\mu(A)^k}{k!}.$$
Any measure-preserving map $T:\left(X,\mathcal{B},\mu\right)\to\left(Y,\mathcal{C},\nu\right)$ %(more generally, any group
%action by measure-preserving transformations)
naturally gives rise to a measure-preserving map $T_{*}:\left(X^{*},\mathcal{B}^{*},\mu^{*}\right)\to\left(Y^{*},\mathcal{C}^{*},\nu^{*}\right)$ by $T_*\gamma = \gamma \circ T^{-1}$. If $T$ is an endomorphism, the dynamical system
$(X^*,\mathcal{B}^*,\mu^*,T_*)$ is the Poisson suspension of
$(X,\mathcal{B},\mu,T)$.

Following \cite{roy_thesis}, the \emph{Poisson entropy} of an infinite-measure-preserving
transformation is defined as the Kolmogorov entropy of the Poisson suspension. This definition gives rise to a new approach for the entropy theory of infinite-measure-preserving transformations. It retains basic properties of Kolmogorov entropy of finite-measure-preserving transformation: If $S$ is a factor of $T$, its Poisson entropy is less than the Poisson entropy of $T$ (Poisson entropy is thus invariant under weak isomorphism), Poisson entropy of $T^n$ is $\left|n\right|$ times Poisson entropy of $T$. The definition of Poisson entropy generalizes to infinite-measure-preserving amenable group actions.

As proved in \cite{roy_thesis}, the Poisson entropy of a probability-preserving transformation is equal to its Kolmogorov entropy. Theorem~\ref{thm:poisson_eq_krengel_quasi_finite} of this paper generalizes this fact: For any quasi-finite transformation, the Poisson entropy is equal to Parry's entropy and Krengel's entropy (this holds in particular for finite-measure-preserving systems).

We recall that if $\left(X,\mathcal{B},\mu,T\right)$ is conservative,
there exists a unique partition of $X$ into $T$-invariant sets $X_{1}$
and $X_{\infty}$, which are the measurable union of finite (resp.
infinite) ergodic components of $\mu$. If $\mu\left(X_{1}\right)=0$
then $T$ is said of type $\mathbf{II}_{\infty}$ and if $\mu\left(X_{\infty}\right)=0$,
of type $\mathbf{II}_{1}$. Only $\mathbf{II}_{\infty}$ systems will
be of interest for us since the $\mathbf{II}_{1}$ case reduces to the
finite measure case. Moreover the possibility to be confronted to
periodic behavior inside the $\mathbf{II}_{1}$ part brings annoying
and uninteresting technical difficulties.

A \emph{factor} of $T$ is a $\sigma$-finite sub-$\sigma$-algebra $\mathcal{F}$ satisfying $T^{-1}\mathcal{F}=\mathcal{F}$. Observe that the  trivial $\sigma$-algebra is never a factor of a $\mathbf{II}_{\infty}$-system.
Remark also that, if $T$ is of type $\mathbf{II}_{\infty}$, then:
\begin{itemize}
\item $\mu$ is continuous;
\item any factor of $T$ is of type $\mathbf{II}_{\infty}$;
\item any $\sigma$-finite sub-$\sigma$-algebra $\mathcal{A}$ satisfying $T^{-1}\mathcal{A}\subset\mathcal{A}$ has no atom;
\item $\left(X^{*},\mathcal{B}^{*},\mu^{*},T_{*}\right)$ is ergodic.
\end{itemize}

We will require some notations and simple results about Poisson measures.
For each $A \in \mathcal{B}$ and $N  \in X^*$, denote by $N(A):X^* \to \NN$ the
random variable on the probability space $(X^*,\mathcal{B}^*,\mu^*)$ which is the (random) measure of the set $A$.
If $A$ has finite measure, $N(A)$ is Poisson distributed with parameter $\mu(A)$.
If $\mu(A)=\infty$, $N(A)=\infty$ $\mu^*$-almost surely.
For a finite or countable partition $\alpha$, we will denote
by $N(\alpha)=(N(A))_{A \in \alpha}$ the random vector of
Poisson random variables corresponding to $\alpha$.
By definition of Poisson suspension, the coordinates of $N(\alpha)$ are independent.
If $\mathcal{C} \subset \mathcal{B}$ is a $\sigma$-algebra, denote by
$\mathcal{C}^*:=\sigma\Bigl(\{N(A):~A \in \mathcal{C}\}\Bigr)$ the sub-$\sigma$-algebra of $\mathcal{B}^*$
generated by the Poisson random variables of $\mathcal{C}$.
For a measurable partition $\alpha$ of $X$, we write $\alpha^*:=(\sigma(\alpha))^*$, sometimes regarding this as a (not necessarily countable) partition of $X^*$.

It is intuitively clear that lack of atoms for a measure implies no ``multiplicities'' in the Poisson space of this measure.
More formally, we have the following standard lemma:
\begin{lemma}\label{lem:small-partitions-separate}
Assuming there are no atoms of positive measure in $(X,\mathcal{B},\mu)$, $\mu^*$-almost surely there are no multiplicities:
$$\mu^*\Bigl( \{ \exists \, x \in X:~ N(\{x\})\ge 2\}\Bigr) = 0 .$$
\end{lemma}

\begin{lemma}\label{lem:star-intersection}
Let $\alpha$, $\beta$ be sub-$\sigma$-algebras of $\mathcal{B}$. Then
\[
\left(\alpha\cap\beta\right)^{*}=\alpha^{*}\cap\beta^{*}\quad\textrm{mod}\;\mu^{*}.\]
\end{lemma}

\begin{proof}
We refer to \cite{roy_poisson} for details about the Fock space structure of $L^2(\mu^*)$ and the exponential $\widetilde\phi$ of an operator $\phi$ of $L^2(\mu)$. For a $\sigma$-algebra $\xi$, let $\pi_{\xi}$ denote the conditional expectation with respect to $\xi$.
It is shown in \cite{roy_poisson} that $\pi_{\alpha^{*}}=\widetilde{\pi_{\alpha}}$
and $\pi_{\beta^{*}}=\widetilde{\pi_{\beta}}$.
Set $H=\left\{ f\in L^{2}\left(\mu\right),\;\pi_{\alpha}\pi_{\beta}f=f\right\} $
and $K=\left\{ g\in L^{2}\left(\mathcal{P}_{\mu}\right),\;\pi_{\alpha^{*}}\pi_{\beta^{*}}g=g\right\} $.
Von Neumann theorem for contractions implies that $\frac{1}{n}{\displaystyle \sum_{k=1}^{n}}\left(\pi_{\alpha}\pi_{\beta}\right)^{k}\to\pi_{H}$
in $L^{2}\left(\mu\right)$ and $\frac{1}{n}{\displaystyle \sum_{k=1}^{n}}\left(\pi_{\alpha^{*}}\pi_{\beta^{*}}\right)^{k}\to\pi_{K}$
in $L^{2}\left(\mathcal{P}_{\mu}\right)$.
But $\pi_{\alpha}\pi_{\beta}f=f$
is equivalent to $\pi_{\alpha}f=\pi_{\beta}f=f$. Therefore
$\pi_{H}=\pi_{\alpha\cap\beta}$.
For the same reason, $\pi_{K}=\pi_{\alpha^{*}\cap\beta^{*}}$.
Moreover, $\widetilde{\frac{1}{n}{\displaystyle \sum_{k=1}^{n}}\left(\pi_{\alpha}\pi_{\beta}\right)^{k}}$
tends to $\widetilde{\pi_{\alpha\cap\beta}}=\pi_{\left(\alpha\cap\beta\right)^{*}}$.
But, for all $n$, $\widetilde{\frac{1}{n}{\displaystyle \sum_{k=1}^{n}}\left(\pi_{\alpha}\pi_{\beta}\right)^{k}}=\frac{1}{n}{\displaystyle \sum_{k=1}^{n}}\left(\pi_{\alpha^{*}}\pi_{\beta^{*}}\right)^{k}$
thus, by uniqueness of the limit, $\pi_{\alpha^{*}\cap\beta^{*}}=\pi_{\left(\alpha\cap\beta\right)^{*}}$,
that is, $\left(\alpha\cap\beta\right)^{*}=\alpha^{*}\cap\beta^{*}$.
\end{proof}

In general, the equality $\left(\mathcal{C}_{1}\vee\mathcal{C}_{2}\right)^{*}=\mathcal{C}_{1}^{*}\vee\mathcal{C}_{2}^{*}$ does not hold. This is however true if the intersection of the $\sigma$-algebras is non-atomic. This is the concern of the next lemma (appearing also in \cite{roy_poisson}):
\begin{lemma}\label{lem:star-partitions}
Let $\alpha$, $\beta$ and $\mathcal{C}$ be sub-$\sigma$-algebras
of $\mathcal{B}$.
Assume that $\mathcal{C}$ is $\sigma$-finite
and non-atomic. Then
\[
\left(\mathcal{C}\vee\alpha\vee\beta\right)^{*}=\left(\mathcal{C}\vee\alpha\right)^{*}\vee\left(\mathcal{C}\vee\beta\right)^{*}\quad\textrm{mod}\;\mu^{*}.\]
\end{lemma}

\begin{proof}
Obviously $\left(\mathcal{C}\vee\alpha\vee\beta\right)^{*}\supset\left(\mathcal{C}\vee\alpha\right)^{*}\vee\left(\mathcal{C}\vee\beta\right)^{*}$.

To complete the proof of this lemma, we must show that for any $A\in\alpha$,
$B\in\beta$ and $C\in\mathcal{C}$, the random variable $N\left(A\cap B\cap C\right)$
is measurable with respect to $\left(\mathcal{C}\vee\alpha\right)^{*}\vee\left(\mathcal{C}\vee\beta\right)^{*}$
up to a $\mu^{*}$-null set. We can find a sequence $(\xi_{n})$ of finite $\mathcal{C}$-measurable partitions increasing to $\mathcal{C}$. Assume $C$ has finite measure. Then by Lemma~\ref{lem:small-partitions-separate}, for almost every $\gamma\in X^{*}$, we consider the smallest integer $n\left(\gamma\right)$
such that for all $E\in\xi_{n\left(\gamma\right)}$, $\gamma\left(E\cap C\right)=0\:\textrm{or}\:1$.
We have
\[
N\left(A\cap B\cap C\right)=\sum_{k\in\mathbb{N}}1_{\left\{ n\left(\gamma\right)=k\right\} }\sum_{E\in\xi_{k}}N\left(A\cap B\cap C\cap E\right).\]

For $E\in\xi_{k}$, set $N^{\prime}\left(E\right):=\min\left(N\left(A\cap C\cap E\right),N\left(B\cap C\cap E\right)\right)$.
Obviously, $N\left(A\cap B\cap C\cap E\right)\leq N^{\prime}\left(E\right)$.
On the other hand, on the set $\left\{ n\left(\gamma\right)=k\right\} $,
$N\left(E\cap C\right)=0\:\textrm{or}\:1$, therefore $N^{\prime}\left(E\right)=1$
if and only if $N\left(A\cap B\cap C\cap E\right)=1$.

Hence, we can write \[
N\left(A\cap B\cap C\right)=\sum_{k\in\mathbb{N}}1_{\left\{ n\left(\gamma\right)=k\right\} }\sum_{E\in\xi_{k}}N^{\prime}\left(E\right).\]

But the right-hand side is measurable with respect to $\left(\mathcal{C}\vee\alpha\right)^{*}\vee\left(\mathcal{C}\vee\beta\right)^{*}$,
so the claim is proved when $C$ has finite measure. In the general case, we can write $C$ as the increasing union of finite-measure, $\mathcal{C}$-measurable sets and get the result in the limit.
\end{proof}

A lemma of the same flavor, which applies to monotone sequences of $\sigma$-algebras, was proved in \cite{roy_thesis}, using the corresponding projections in $L^2(\mu)$ and $L^2(\mu^*)$.
\begin{lemma}\label{lem:star_monotone_algebras_roy}
Let $\{\mathcal{B}_n\}_{n \in \NN}$ be a sequence of sub-$\sigma$-algebras of $\mathcal{B}$.
\begin{enumerate}
\item{If $\{\mathcal{B}_n\}_{n \in \NN}$ is an increasing sequence,
then $\bigvee_{n \in \NN}\mathcal{B}_n^* = (\bigvee_{n \in\NN}\mathcal{B}_n)^*$.}
\item{If $\{\mathcal{B}_n\}_{n \in \NN}$ is a decreasing sequence,
then $\bigcap_{n \in \NN}\mathcal{B}_n^* = (\bigcap_{n \in\NN}\mathcal{B}_n)^*$.}
\end{enumerate}
The above equalities are modulo null sets.
\end{lemma}

\section{The Krengel entropy of a conservative measure-preserving transformation}

The Krengel entropy of a conservative measure-preserving
transformation $(X,\mathcal{B},\mu,T)$ is defined
in~\cite{krengel67} as:
$$
\hkr (X,\mathcal{B},\mu,T) :=\sup_{A \in\mathcal{F}_+}\mu(A)\,h(A,\mathcal{B}\cap A,{\mu}_A,T_A),
$$
where $\mathcal{F}_+$ is the collection of sets in $\mathcal{B}$ with
finite positive measure, ${\mu}_A$ is the normalized
probability measure on $A$ obtained by restricting $\mu$ to
$\mathcal{B} \cap A$, and $T_A:A \to A$ is the induced map on
$A$. Recall that this map is defined by
$$T_A(x) := T^{\phi_A(x)}(x),$$
where $\phi_A(x) := \min\{k \ge 1:~ T^k(x) \in A\}$ is the \emph{first-return-time map} associated to $A$.
As soon as $T$ is not purely periodic, Krengel proved that
$$
\hkr (X,\mathcal{B},\mu,T) =\mu(A)\,h(A,\mathcal{B}\cap A,\mu_A,T_A),
$$
where $A$ is any finite-measure \emph{sweep-out} set (i.e. a set such that
${\displaystyle \bigcup_{n=0}^{\infty}}T^{-n}A=X$), which always exists when $T$ is of type $\mathbf{II}_{\infty}$.

The fact that Krengel's entropy extends Kolmogorov's follows from
Abramov's formula.  The latter states that when $S:\Omega \to \Omega$ is an ergodic probability-preserving transformation on $(\Omega,\mathcal{F},p)$, and $A\in \mathcal{B}$, we have
$$h(A,\mathcal{F} \cap A, p(\cdot \mid A),S_A)=\frac{1}{p(A)}h(\Omega,\mathcal{F},p,S).$$

\section{\label{sec:Info}The information function of a measurable partition}

We describe here a generalization of Shannon's information function. This was previously studied by Klimko and Sucheston \cite{klimko_sucheston68} and Parry \cite{parry_generators69}:
The \emph{information function} of a partition $\alpha$ is given by
$$
I_{\mu}(\alpha)(x) := \begin{cases}
\log\frac{1}{\mu(\alpha(x))} & \mbox{if } 0 < \mu(\alpha(x)) < \infty,\\
\infty & \mbox{if } \mu(\alpha(x))=0, \\
0 &  \mbox{if } \mu(\alpha(x))=\infty .
\end{cases}
$$
By $\alpha(x)$ we mean the unique element in $\alpha$ which contains $x$.
Similarly, given two partitions $\alpha_1$ and $\alpha_2$, the \emph{conditional information} is defined as
\begin{equation}
\label{eq:information_func}
I_{\mu}(\alpha_1 \mid \alpha_2)(x) :=
\begin{cases}
I_{\mu(\cdot \mid \alpha_2(x))}(\alpha_1 )(x) & \mbox{if }\mu(\alpha_2(x)< \infty),\\
I_{\mu}\Bigl(\alpha_1 \vee \{(\alpha_2(x)),X\setminus \alpha_2(x)\}\Bigr) (x)& \mbox{otherwise.}
\end{cases}
\end{equation}

Note that the conditional information retains the following property from the finite-measure case (see \cite{klimko_sucheston68}):
\begin{equation}
 \label{information_decomposition}
I_{\mu}\left(\bigvee_{0}^{n-1}T^k\alpha\right)
= I_{\mu}\left(\alpha\right)\circ T^{-\left(n-1\right)}+\sum_{j=1}^{n-1}I_{\mu}\left(\alpha \left| \bigvee_{1}^{j}T^k
\right.\right)\circ T^{j-\left(n-1\right)}
\end{equation}

In the sequel, we will need the following lemma (see Theorem~2.2 in \cite{parry_generators69} for a proof):
\begin{lemma}\label{lem:info_mean_converge}
Let $(\Omega,\mathcal{F},p)$ be a probability space, $\alpha$ a measurable partition with $H_{\mu}(\alpha) <\infty$, and $\{\mathcal{F}_n\}_{n \in \mathbb{N}}$ an increasing sequence of sub-$\sigma$-algebras such that $\mathcal{F}=\bigvee_{n \ge 1}\mathcal{F}_n$. Then
$$I_\mu(\alpha \mid \mathcal{F}_n) \to I_\mu(\alpha \mid \mathcal{F})$$ in $L_1(\Omega,p)$ and $p$-a.e.
\end{lemma}

When $\mathcal{C}$ and $\mathcal{D}$ are sub-$\sigma$-algebras
corresponding to partitions $\alpha$ and $\beta$, we note
$I_{\mu}(\mathcal{C}) := I_{\mu}( \alpha)$ and
$I_{\mu}(\mathcal{C} \mid \mathcal{D}) := I_{\mu}(\alpha \mid \beta)$.

Following Parry, if $\mathcal{C}$ and $\mathcal{D}$ are
$\sigma$-finite sub-$\sigma$-algebras of $\mathcal{B}$, we define the entropy of $\mathcal{C}$ by
$$H_{\mu}(\mathcal{C}) :={\displaystyle \int_{X}}I_{\mu}\left(\mathcal{C}\right)d\mu$$ and the
\emph{conditional entropy} of $\mathcal{C}$ with respect to $\mathcal{D}$ by
$$H_{\mu}(\mathcal{C} \mid \mathcal{D})
:= \int_X I_\mu(\mathcal{C} \mid \mathcal{D})(x) d\mu(x).
$$
Since $\mathcal{D}$ is $\sigma$-finite, we have $\mu(\beta(x)< \infty)$ for $\mu$-almost all $x$ where $\beta$ is the partition associated to $\mathcal{D}$. Hence,
$$H_{\mu}(\mathcal{C} \mid \mathcal{D}) = \int_X H_{\mu(\cdot|\beta(x))}(\mathcal{C})\, d\mu(x).$$

Finite conditional entropy implies that $\mu$-almost every atom of
$\mathcal{D}$ intersects at most countably many atoms of $\mathcal{C}$.

The following lemma is useful for entropy estimates of a Poisson measure.
\begin{lemma}\label{lem:poisson_cond_entropy}
Assume that $\left(X,\mathcal{B},\mu\right)$ is a Lebesgue space where $\mu$ is continuous and infinite. Let $\mathcal{D} \subset \mathcal{C}$ be
$\sigma$-finite sub-$\sigma$-algebras of $\mathcal{B}$ with no atom of positive
measure. Then
$$H_{\mu^*}(\mathcal{C}^* \mid \mathcal{D}^*) = H_{\mu}(\mathcal{C} \mid \mathcal{D}).$$
\end{lemma}
\begin{proof}
Note first that we can take $\mathcal{C}=\mathcal{B}$ and by disintegrating
$\mu$ with respect to $\mathcal{D}$ and using the fact that $\left(X,\mathcal{B},\mu\right)$
is a Lebesgue space with a continuous infinite measure, we can represent $\left(X,\mathcal{B},\mu\right)$
as $\left(\mathbb{R}\times Y,\mathcal{A}\otimes\mathcal{Y},\mu\right)$,
where $\mathcal{A}$ is the Borel $\sigma$-algebra and
$\mu\left(A_{1}\times A_{2}\right)={\displaystyle \int_{A_{1}}}m_{x}\left(A_{2}\right)\lambda\left(dx\right)$,
$\lambda$ being the Lebesgue measure.
Thus $\left(X^{*},\mathcal{B}^{*},\mu^{*}\right)$ takes the form
$\left(\left(\mathbb{R}\times Y\right)^{*},\left(\mathcal{A}\otimes\mathcal{Y}\right)^{*},\mu^{*}\right)$.
This latter Poisson measure has the form of a so-called \emph{marked Poisson process}, namely, we can identify it with
$\left(\mathbb{R}^{*}\times Y^{\mathbb{Z}},\mathcal{A}^{*}\otimes\mathcal{Y}^{\otimes\mathbb{Z}},\mathbb{P}\right)$
through the mapping $\nu\in\left(\mathbb{R}\times Y\right)^{*}\mapsto\left(\gamma,\left\{ y_{i}\right\} _{i\in\mathbb{Z}}\right)$,
where
\begin{itemize}
 \item $\gamma$ is the projection of $\nu$ on $\mathbb{R}$:
$\gamma= {\displaystyle \sum_{i\in\mathbb{Z}}}\delta_{t_{i}\left(\gamma\right)}$, with $$\cdots<t_{-n}\left(\gamma\right)<\cdots<t_{-1}\left(\gamma\right)<t_{0}\left(\gamma\right)\leq0<t_{1}\left(\gamma\right)<\cdots<t_{n}\left(\gamma\right)<\cdots$$
\item $(y_{i})_{i\in\mathbb{Z}}$ are defined by $\nu={\displaystyle \sum_{i\in\mathbb{Z}}}\delta_{\left(t_{i}\left(\gamma\right),y_{i}\right)}$.
\end{itemize}
Here $\mathbb{P}\left(C_{1}\times C_{2}\right)={\displaystyle \int_{C_{1}}}p_{\gamma}\left(C_{2}\right)\lambda^{*}\left(d\gamma\right)$
with $p_{\gamma}=\otimes_{i\in\mathbb{Z}}m_{t_{i}\left(\gamma\right)}$.
The verification of this fact amounts to evaluating the Laplace transform for a positive function $f$ on $\mathbb{R}\times{Y}$.
%$\int_{\mathbb{R}^{*}\times Y^{\mathbb{Z}}}\exp\left(- \sum_{i\in\mathbb{Z}} f\left(t_{i}(\gamma),y_{i}\right) \right)\mathbb{P}\left(d\left(\gamma,\left\{ y_{i}\right\} _{i\in\mathbb{Z}}\right)\right)$
\begin{eqnarray*}
\lefteqn{\int_{\mathbb{R}^{*}}\left( \int_{Y^{\mathbb{Z}}} \exp\left(- \sum_{i\in\mathbb{Z}} f\left(t_{i}(\gamma),y_{i}\right)\right)\otimes_{i\in\mathbb{Z}}m_{t_{i}\left(\gamma\right)}\left(d\left\{ y_{i}\right\} _{i\in\mathbb{Z}}\right)\right)\lambda^{*}\left(d\gamma\right)}\\
%&=& \int_{\mathbb{R}^{*}} \prod_{i\in\mathbb{Z}}\left( \int_{Y} \exp\left(-f\left(t_{i}(\gamma),y\right)\right)m_{t_{i}\left(\gamma\right)}\left(dy\right)\right)\mu^{*}\left(d\gamma\right)\\
&=& \int_{\mathbb{R}^{*}}\exp \left\{\sum_{i\in\mathbb{Z}}\log\left( \int_{Y} \exp\left(-f\left(t_{i}(\gamma),y\right)\right)m_{t_{i}\left(\gamma\right)}\left(dy\right)\right)\right\}\lambda^{*}\left(d\gamma\right)\\
&=& \int_{\mathbb{R}^{*}}\exp \left\{\int_{\mathbb{R}} \log\left( \int_{Y} \exp\left(-f(t,y)\right)m_{t} \left(dy\right)\right)d\gamma\left(t\right)\right\}\lambda^{*}\left(d\gamma\right)\\
% &=& \int_{\mathbb{R}^{*}}\left(\exp \int_{\mathbb{R}} \left(\log\left( \int_{Y} \exp\left(-f\left(t,y\right)\right)m_{t} \left(dy\right)\right)\right)d\gamma\left(t\right)\right)\lambda^{*}\left(d\gamma\right)\\
&=&\exp \int_{\mathbb{R}} \left\{ \exp\left(\log\left( \int_{Y} \exp\left(-f\left(t,y\right)\right)m_{t}\left(dy\right)\right)\right)-1 \right\} dt\\
%&=&\exp \int_{\mathbb{R}} \left( \int_{Y} \exp\left(-f\left(t,y\right)\right)m_{t}\left(dy\right)\right)-1dt\\
&=&\exp \int_{\mathbb{R}} \left( \int_{Y} \left(\exp\left(-f\left(t,y\right)\right)-1\right)m_{t}\left(dy\right)\right)dt\\
&=&\exp \int_{\mathbb{R}\times Y} \left(\exp\left(-f\left(z\right)\right)-1\right)\mu\left(dz\right)\\
&=&\int_{\left(\mathbb{R}\times Y\right)^{*}}\exp\left(-\int_{\mathbb{R}\times Y}f\left(z\right)\rho\left(dz\right)\right)\mu^{*}\left(d\rho\right),
\end{eqnarray*}
which is the Laplace transform of the Poisson measure of distribution $\mu^{*}$ evaluated at $f$.

% The verification amounts to evaluating $\mathbb{P}\left(\left\{\left(\gamma,\left\{ y_{i}\right\} _{i\in\mathbb{Z}}\right): {\displaystyle \sum_{i\in\mathbb{Z}}}\delta_{\left(t_{i}\left(\gamma\right),y_{i}\right)}\left(\left[a,b\right]\times C\right)=0\right\}\right)$
% for any interval $\left[a,b\right]\subset\mathbb{R}$ and any $C\in\mathcal{Y}$. \textbf{justifications ???}

In this setting, we can rewrite $H_{\mu^{*}}\left(\mathcal{C}^{*}\mid\mathcal{D}^{*}\right)$ as
$$
\int_{\mathbb{R}^{*}} d\lambda^*(\gamma)\, H_{p_{\gamma}}\left( \mathcal{Y}^{\otimes\mathbb{Z}}\right)
= \int_{\mathbb{R}^{*}} d\lambda^*(\gamma)\,  \sum_{i\in\mathbb{Z}} H_{m_{t_i}}(\mathcal{Y})
= \int_{\mathbb{R}^{*}}  d\lambda^*(\gamma)\, \int_{\mathbb{R}}d\gamma(t)\, H_{m_{t}}(\mathcal{Y}),
$$
which is equal to
$$
\int_{\mathbb{R}} d\lambda(t)\, H_{m_{t}}(\mathcal{Y})
=H_{\mu}\left(\mathcal{C}\mid\mathcal{D}\right).
$$

\end{proof}

\section{\label{sec:Parry} Parry's entropy}
In this section we recall Parry's definition of entropy for a measure-preserving transformation, and prove that Parry's entropy is dominated by Poisson entropy.

Parry \cite{parry_generators69} defines the entropy of a measure-preserving transformation by
$$
\hpa(X,\mathcal{B},\mu,T) := \sup_{T^{-1}\mathcal{C} \subset \mathcal{C}}H_{\mu}(\mathcal{C} \mid T^{-1}\mathcal{C}),
$$
where the supremum is taken over all $\sigma$-finite sub-$\sigma$-algebras $\mathcal{C}$
of $\mathcal{B}$ such that $T^{-1}\mathcal{C} \subset \mathcal{C}$.
For probability-preserving transformations, this definition coincides with the standard definition of Kolmogorov's entropy.

The following theorem was proved by Parry (Theorem~10.11 in \cite{parry_generators69}).
\begin{theorem}
\label{thm:poisson_le_krnegel}%^^^
Let $(X,\mathcal{B},\mu,T)$ be a measure-preserving conservative transformation. Then
$$ \hpa(X,\mathcal{B},\mu,T) \le \hkr(X,\mathcal{B},\mu,T). $$
%the Parry entropy of $(X,\mathcal{B},\mu,T)$ is dominated by the Krengel entropy of $(X,\mathcal{B},\mu,T)$.
\end{theorem}

Replacing Krengel entropy by Poisson entropy, we prove a similar result:
\begin{theorem}\label{thm:poisson_eq_parry}
Let $(X,\mathcal{B},\mu,T)$ be a  $\mathbf{II}_{\infty}$ transformation.
Then
$$ \hpa(X,\mathcal{B},\mu,T) \le h(X^*,\mathcal{B}^*,\mu^*,T_*). $$
%the Parry entropy of $(X,\mathcal{B},\mu,T)$ is dominated by the Poisson entropy.
\end{theorem}
\begin{proof}
Let $\mathcal{C} \subset \mathcal{B}$ be a sub-invariant
$\sigma$-finite sub-$\sigma$-algebra, that is $T^{-1}\mathcal{C} \subset \mathcal{C}$.
Since $\mathcal{B}$ has no atom of positive $\mu$-measure and $\mathcal{C}$ is $\sigma$-finite,
the same follows for $\mathcal{C}$, and so by Lemma~\ref{lem:poisson_cond_entropy} we
know that
$$
H_{\mu}(\mathcal{C} \mid T^{-1}\mathcal{C}) = H_{\mu^*}(\mathcal{C}^* \mid T_*^{-1}\mathcal{C}^*).
$$
Now it follows that
$$\sup_{T^{-1}\mathcal{C} \subset \mathcal{C}}H_{\mu}(\mathcal{C} \mid T^{-1}\mathcal{C}) \le
\sup_{T_*^{-1}\mathcal{D} \subset \mathcal{D}}H_{\mu^*}(\mathcal{D} \mid T^{-1}\mathcal{D}), $$
where the supremum on the right-hand side is over all factors $\mathcal{D} \subset \mathcal{B}^*$, which proves the theorem.
\end{proof}

\section{An upper bound for the Poisson entropy}

Whenever the measure-preserving system $(X,\mathcal{B},\mu,T)$ is implicitly clear from the context,
for any measurable partition $\alpha$ of $X$ and $-\infty \le i<j \le +\infty$, we write $\alpha^j_i := \bigvee_{k=i}^j T^{k}\alpha$.
We will assume from now on that $T$ is an automorphism, that is $T^{-1}\mathcal{B}= \mathcal{B}$ with equality modulo $\mu$.
Also, we write $\widehat{\alpha}=\alpha_{-\infty}^{\infty}$.

We say that a countable partition $\alpha$ of $X$ is
\emph{local} with \emph{core} $A \in \mathcal{F}$ if
$$A^c \in \alpha \quad\mbox{ and } \quad  H_{\mu}(\alpha) < \infty\, .$$
In other words, $\alpha$ is a finite-entropy partition of a set $A$ of finite measure, to which the complement of $A$ is added.
Note that $\alpha^*$ is at most a countable partition of $X^*$, since %$H_{\mu}(\alpha)< \infty$ implies $H_{\mu^*}(\alpha^*)<\infty$, and a partition with finite entropy is always (essentially) countable.
$N(A)$ is finite $\mu^*$-almost surely.
In \cite{aaro_park} and \cite{klimko_sucheston68} \emph{co-finite} partitions were considered: With our terminology, these are finite, local partitions.

For a  partition $\alpha$ of $X$, define
$$\hat{h}(X,\mu,T,\alpha) := \liminf_{n \to\infty}\frac{1}{n}H_{\mu}(\alpha_0^{n-1}).$$
In case $\mu$ is a probability measure and $H_{\mu}(\alpha)<\infty$, this is equal to the Kolmogorov entropy of the factor generated by $\{T^{-n}\alpha\}_{n=0}^{\infty}$.

% We also define, for any measurable partition $\alpha$,
% $$\hat{\lambda}(\alpha) := \sup_{a \in \alpha \cap \mathcal{F}}\mu(a).$$

\begin{lemma}\label{lem:hat_partition_parameter}
Let $(X,\mathcal{B},\mu,T)$ be a $\mathbf{II}_{\infty}$ transformation
and let $\alpha$ be a local partition whose core $A$ is a sweep-out set.
Then we have
$$
\lim_{n \to \infty} \sup \left\{ \mu \left(a\right) : a \in  \bigvee_{k=0}^{n-1}T^{-k}\alpha\cap \mathcal{F} \right\} = 0 .
$$
\end{lemma}
\begin{proof}
By considering the natural extension of $(X,\mathcal{B},\mu,T)$, we can assume that the transformation is invertible.
Since $A$ is a sweep-out set, the first-return-time map $\phi_A$ is finite almost everywhere.
We claim that the first-return-time map $\psi_A(x)=\min\{n >0:~T^{-n}(x) \in A\}$ for $T^{-1}$ is also finite almost everywhere.
Indeed, let
$$ C := \{x\in X : \ \forall n>0, T^{-n}(x) \notin A\}. $$
For all $x\in X$, the number of positive $n$'s such that $T^nx\in C$ is bounded by $\phi_A(x)$, thus is finite almost everywhere. Since $T$ is conservative, this implies that $\mu(C)=0$.

Choose $k \in \NN$ so large that $\mu(A \cap \{\phi_A(x) >k\}) < \epsilon$ and
$\mu(A \cap \{\psi_A(x) >k\}) < \epsilon$.
For all $n\ge1$, let
$$B_{n}:=B\cap T_A^{-1}B\cap\cdots\cap T_A^{-n}B,$$
where
$$B:=\{x \in A:~ \phi_A(x) < k \}.$$
We claim that $\lim_{n\to \infty}\mu(B_n)=0$. Indeed, if $\lim_{n\to \infty}\mu(B_n)>0$, then
$B_{\infty}:=\bigcap_{n\ge1}B_{n}$ is a set of positive measure, $T_A$-invariant, and its first-return-time map $\phi_{B_\infty}$ is bounded by $k$. Then $B_{\infty}\cup TB_{\infty}\cup\cdots \cup T^kB_{\infty}$ is a set of finite positive measure which is $T$-invariant, which contradicts the hypothesis that $(X,\mathcal{B},\mu,T)$ is $\mathbf{II}_{\infty}$.

We conclude that every $a \in \bigvee_{j=0}^{kn}T^{-j}\mathcal{\alpha} \cap \mathcal{F}$ is either
contained in $T^{-j}[A \cap\{\phi_A >k\}]$ or in $T^{-j}[A \cap \{\psi_A >k \}]$ for some $j \in \NN$, or $a \subset B_{n}$.
If $n$ is large enough, we get that $\mu(a)<\epsilon$.
\end{proof}

\begin{proposition}\label{prop:hat_ge_poisson}
Let $(X,\mathcal{B},\mu,T)$ be a $\mathbf{II}_{\infty}$ transformation
and let $\alpha$ be a local partition, whose core is a sweep-out set. We have
$$h(X^*,(\hat{\alpha})^*,\mu^*,T_*)\le \hat{h}(X,\mu,T,\alpha). $$
\end{proposition}
\begin{proof}
For any $n\geq1$ and any $p\geq1$, since $\left(\left(\alpha_{0}^{p}\right)^{*}\right)_{0}^{n-1}\prec\left(\alpha_{0}^{n-1+p}\right)^{*}$, we have
$$
\frac{1}{n} H_{\mu^{*}}\left(\left(\left(\alpha_{0}^{p}\right)^{*}\right)_{0}^{n-1}\right)
\leq \frac{n+p}{n} \frac{1}{n+p} H_{\mu^{*}}\left(\left(\alpha_{0}^{n-1+p}\right)^{*}\right).
$$
But $H_{\mu^{*}}\left(\left(\alpha_{0}^{n-1+p}\right)^{*}\right)=\sum_{a\in\alpha_{0}^{n-1+p}}f\left(\mu\left(a\right)\right)$
where $f\left(x\right)$ is the entropy of a Poisson random variable
with parameter $x$.
An easy computation shows that $f\left(\epsilon\right)\sim-\epsilon\log\epsilon$
at the origin.
By Lemma~\ref{lem:hat_partition_parameter}, %$\hat{\lambda}\left(\alpha_{0}^{n-1+p}\right)$
$\sup_{a \in \alpha_{0}^{n-1+p} \cap \mathcal{F}}\,\mu(a)$ tends to $0$ as $n$ tends to infinity. Therefore, $H_{\mu^{*}}\left(\left(\alpha_{0}^{n-1+p}\right)^{*}\right) \sim H_{\mu}\left(\alpha_{0}^{n-1+p}\right)$
as $n$ tends to infinity, and we get
\[
\lim_{n\to\infty} \frac{1}{n} H_{\mu^{*}}\left(\left(\left(\alpha_{0}^{p}\right)^{*}\right)_{0}^{n-1}\right)
\leq \liminf_{n\to\infty} \frac{n+p}{n} \frac{1}{n+p} H_{\mu}\left(\alpha_{0}^{n-1+p}\right)
=\hat{h}\left(X,\mu,T,\alpha\right).\]
Taking the limit in $p$, we obtain the desired inequality.
\end{proof}

\section{Relative Poisson entropy}

Here $\mathcal{C}$ is an invertible factor  ($T^{-1}\mathcal{C}=\mathcal{C}$) of $\left(X,\mathcal{B},\mu,T\right)$.
The \emph{relative entropy of $T$ with respect to $\mathcal{C}$} is defined by % in \cite{danilanko_rudolph07}:
\begin{equation}\label{eq:relative_entropy}
\hrel(X,\mathcal{B},\mu,T \mid \mathcal{C}) :=
\sup_{\alpha}\lim_{n \to \infty} \frac{1}{n} H_{\mu}\left( \bigvee_{k=0}^{n-1}T^{-k}\alpha \Big\vert  \mathcal{C} \right) ,
\end{equation}
where the supremum is taken over all countable partitions $\alpha$ with
$H_{\mu}(\alpha \mid \mathcal{C}) < \infty$. %, or some dense subset of them.

This definition of relative entropy, which is classical for
probability-preserving transformations, was applied to
$\sigma$-finite measure-preserving actions of countable amenable groups by
 Danilenko and Rudolph \cite{danilenko_rudolph07}.

\begin{proposition}\label{prop:relative_entropy}
Let $\mathcal{C}$ be an invertible factor of a $\mathbf{II}_{\infty}$ system $\left(X,\mathcal{B},\mu,T\right)$.
Then the following quantities are equal:
\begin{itemize}
\item $\hrel\left(X,\mathcal{B},\mu,T\mid\mathcal{C}\right)$
\item ${\displaystyle \lim_{p\to\infty}}{\displaystyle \lim_{n\to\infty}}\frac{1}{n}H_{\mu}\left({\displaystyle \bigvee_{k=0}^{n-1}}T^{-k}\alpha_{p}\Big|\mathcal{C}\right)$, where $\alpha_{p}\uparrow\mathcal{B}$  are a sequence of local partitions with a core $A \in \mathcal{C}$ which is a sweep out set
    and
$H_{\mu}\left(\alpha_{p}\mid\mathcal{C}\right)<\infty$.
\item ${\displaystyle \sup_{\mathcal{D}\subset\mathcal{B},T^{-1}\mathcal{D}\subset\mathcal{D}}}H_{\mu}\left(\mathcal{D}\mid T^{-1}\mathcal{D}\vee\mathcal{C}\right)$ ($\mathcal{D}$ $\sigma$-finite)
\item ${\displaystyle \sup_{\mathcal{D}\subset\mathcal{B},T_{*}^{-1}\mathcal{D}^{*}\subset\mathcal{D}^{*}}}H_{\mu}\left(\mathcal{D}^{*}\mid T_{*}^{-1}\mathcal{D}^{*}\vee\mathcal{C}^{*}\right)$ ($\mathcal{D}$ $\sigma$-finite)
\item $h\left(X^{*},\mathcal{B}^{*},\mu^{*},T_{*}\mid\mathcal{C}^{*}\right)$
\item $\mu\left(A\right)h\left(A,\mathcal{B}\cap A,\mu_{A},T_{A}\mid\mathcal{C}\cap A\right)$
for any sweep-out set $A\in\mathcal{C}$.
\end{itemize}
\end{proposition}
\begin{proof}
Let $\alpha$ be a local partition whose core $A\in\mathcal{C}$ is a sweep out set, and such that $H_{\mu}\left(\alpha\mid\mathcal{C}\right)<\infty$.

\[
H_{\mu}\left({\displaystyle \bigvee_{k=0}^{n-1}}T^{-k}\alpha\Big|\mathcal{C}\right)=\sum_{k=0}^{n-1}H_{\mu}\left(T^{-k}\alpha\Big|{\displaystyle \bigvee_{j=k+1}^{n-1}}T^{-j}\alpha\vee\mathcal{C}\right)\]

%
% \[
% =\sum_{k=0}^{n-1}H_{\mu}\left(T^{-k}\alpha\Big| T^{-k}\left({\displaystyle \bigvee_{j=1}^{n-1-k}}T^{-j}\alpha\right)\vee T^{-k}\mathcal{C}\right)\]

\[
=\sum_{k=0}^{n-1}H_{\mu}\left(\alpha\Big|{\displaystyle \bigvee_{j=1}^{n-1-k}}T^{-j}\alpha\vee\mathcal{C}\right)
=\sum_{k=0}^{n-1}H_{\mu}\left(\alpha\Big|{\displaystyle \bigvee_{j=1}^{k}}T^{-j}\alpha\vee\mathcal{C}\right)
\]
Since $H_{\mu}\left(\alpha\Big|{\displaystyle \bigvee_{j=1}^{k}}T^{-j}\alpha\vee\mathcal{C}\right)$
tends to \[
H_{\mu}\left(\alpha\Big|{\displaystyle \bigvee_{j=1}^{\infty}}T^{-j}\alpha\vee\mathcal{C}\right)=H_{\mu}\left({\displaystyle \bigvee_{j=0}^{\infty}}T^{-j}\alpha\Big|{\displaystyle \bigvee_{j=1}^{\infty}}T^{-j}\alpha\vee\mathcal{C}\right),\]
Cesaro averages gives
\[
{\displaystyle \lim_{n\to\infty}}\frac{1}{n}H_{\mu}\left({\displaystyle \bigvee_{k=0}^{n-1}}T^{-k}\alpha\Big|\mathcal{C}\right)=H_{\mu}\left({\displaystyle \bigvee_{j=0}^{\infty}}T^{-j}\alpha\Big|{\displaystyle \bigvee_{j=1}^{\infty}}T^{-j}\alpha\vee\mathcal{C}\right).
\]
Now remark that, since $A\in\mathcal{C}$, $I_{\mu}\left(\alpha\Big|{\displaystyle \bigvee_{j=1}^{\infty}}T^{-j}\alpha\vee\mathcal{C}\right)$
vanishes outside $A$ and $\left({\displaystyle \bigvee_{j=1}^{\infty}}T^{-j}\alpha\vee\mathcal{C}\right)\cap A=\left({\displaystyle \bigvee_{j=1}^{\infty}}T_{A}^{-j}\alpha\vee\mathcal{C}\right)\cap A$.
Thus

\begin{eqnarray*}
 &  & H_{\mu}\left(\alpha\Big|{\displaystyle \bigvee_{j=1}^{\infty}}T^{-j}\alpha\vee\mathcal{C}\right)={\displaystyle \int_{X}}I_{\mu}\left(\alpha\Big|{\displaystyle \bigvee_{j=1}^{\infty}}T^{-j}\alpha\vee\mathcal{C}\right)d\mu\\
 &  & =\mu(A){\displaystyle \int_{A}}I_{\mu_{A}}\left(\alpha\Big|\left({\displaystyle \bigvee_{j=1}^{\infty}}T_{A}^{-j}\alpha\vee\mathcal{C}\right)\cap A\right)d\mu_{A}=\mu(A)H_{\mu_{A}}\left(\alpha\Big|\left({\displaystyle \bigvee_{j=1}^{\infty}}T_{A}^{-j}\alpha\vee\mathcal{C}\right)\cap A\right)
\end{eqnarray*}
On the one hand, ${\displaystyle \sup_{\alpha}}H_{\mu_{A}}\left(\alpha\Big|\left({\displaystyle \bigvee_{j=1}^{\infty}}T_{A}^{-j}\alpha\vee\mathcal{C}\right)\cap A\right)$
over countable partition of $A$ such that $H_{\mu}\left(\alpha\mid\mathcal{C}\right)=H_{\mu_{A}}\left(\alpha\mid\mathcal{C}\cap A\right)<\infty$
equals $h\left(A,\mathcal{B}\cap A,\mu_{A},T_{A}\mid\mathcal{C}\cap A\right)$.
On the other hand, $H_{\mu}\left(\alpha\Big|{\displaystyle \bigvee_{j=1}^{\infty}}T^{-j}\alpha\vee\mathcal{C}\right)\leq{\displaystyle \sup_{\mathcal{D}\subset\mathcal{B},T^{-1}\mathcal{D}\subset\mathcal{D}}}H_{\mu}\left(\mathcal{D}\mid T^{-1}\mathcal{D}\vee\mathcal{C}\right)$.
Therefore,
\[
\mu\left(A\right)h\left(A,\mathcal{B}\cap A,\mu_{A},T_{A}\mid\mathcal{C}\cap A\right)\leq \hrel\left(X,\mathcal{B},\mu,T\mid\mathcal{C}\right)\leq{\displaystyle \sup_{\mathcal{D}\subset\mathcal{B},T^{-1}\mathcal{D}\subset\mathcal{D}}}H_{\mu}\left(\mathcal{D}\mid T^{-1}\mathcal{D}\vee\mathcal{C}\right)
\]
Now, since $T$ is of type  $\mathbf{II}_{\infty}$, observe that if $\mathcal{D}$ is a sub-invariant $\sigma$-finite sub-$\sigma$-algebra of $\mathcal{B}$, then $T^{-1}\mathcal{D}$ is non-atomic. We thus have
\begin{eqnarray*}
\lefteqn{H_{\mu}\left(\mathcal{D}\mid T^{-1}\mathcal{D}\vee\pi^{-1}\mathcal{C}\right)} \\
& = & H_{\mu}\left(\mathcal{D}\vee\pi^{-1}\mathcal{C}\mid T^{-1}\mathcal{D}\vee\pi^{-1}\mathcal{C}\right) \\
& = & H_{\mu^*}\left( \left(\mathcal{D}\vee\pi^{-1}\mathcal{C}\right)^{*} \mid \left(T^{-1}\mathcal{D}\vee\pi^{-1}\mathcal{C}\right)^* \right) \quad\mbox{by Lemma~\ref{lem:poisson_cond_entropy}} \\
& = & H_{\mu^{*}}\left(\left( \mathcal{D} \vee T^{-1}\mathcal{D} \vee \pi^{-1}\mathcal{C} \right)^{*} \mid\left(T^{-1}\mathcal{D}\vee\pi^{-1}\mathcal{C}\right)^{*}\right) \\
& = & H_{\mu^{*}}\left(\mathcal{D}^{*}\vee\left(T^{-1}\mathcal{D}\vee\pi^{-1}\mathcal{C}\right)^{*}\mid\left(T^{-1}\mathcal{D} \vee \pi^{-1}\mathcal{C}\right)^{*}\right) \quad\mbox{by Lemma~\ref{lem:star-partitions}} \\
& = &H_{\mu^{*}}\left(\mathcal{D}^{*}\mid\left(T^{-1}\mathcal{D}\vee\pi^{-1}\mathcal{C}\right)^{*}\right) \\
&\le&H_{\mu^{*}}\left(\mathcal{D}^{*}\mid T_{*}^{-1}\mathcal{D}^{*}\vee\pi_{*}^{-1}\mathcal{C}^{*}\right).
\end{eqnarray*}
Hence,
\[
{\displaystyle \sup_{\mathcal{D}\subset\mathcal{B},T^{-1}\mathcal{D}\subset\mathcal{D}}}H_{\mu}\left(\mathcal{D}\mid T^{-1}\mathcal{D}\vee\mathcal{C}\right)\]

\[
\leq{\displaystyle \sup_{\mathcal{D}\subset\mathcal{B},T_{*}^{-1}\mathcal{D}^{*}\subset\mathcal{D}^{*}}}H_{\mu}\left(\mathcal{D}^{*}\mid T_{*}^{-1}\mathcal{D}^{*}\vee\mathcal{C}^{*}\right)\leq h\left(X^{*},\mathcal{B}^{*},\mu^{*},T_{*}\mid\mathcal{C}^{*}\right)\]
Moreover, by taking an increasing sequence $\alpha_{p}$
of countable partitions with core $A$, $H_{\mu}\left(\alpha_{p}\mid\mathcal{C}\right)<\infty$,
such that $\alpha_{p}\cap A\uparrow\mathcal{B}\cap A$, we have $\widehat{\alpha_{p}^{*}}\uparrow\mathcal{B}^{*}$
and therefore
\[
h\left(X^{*},\mathcal{B}^{*},\mu^{*},T_{*}\mid\mathcal{C}^{*}\right)={\displaystyle \lim_{p\to\infty}}{\displaystyle \lim_{n\to\infty}}\frac{1}{n}H_{\mu^{*}}\left({\displaystyle \bigvee_{k=0}^{n-1}}T_{*}^{-k}\alpha_{p}^{*}\Big|\mathcal{C}^{*}\right).\]
But
\begin{eqnarray*}
 &  & \frac{1}{n}H_{\mu^{*}}\left({\displaystyle \bigvee_{k=0}^{n-1}}T_{*}^{-k}\alpha_{p}^{*}\Big|\mathcal{C}^{*}\right)\leq \frac{1}{n}H_{\mu^{*}}\left(\left({\displaystyle \bigvee_{k=0}^{n-1}}T^{-k}\alpha_{p}\right)^{*}\vee\mathcal{C}^{*}\Big|\mathcal{C}^{*}\right)\\
 &  & \leq \frac{1}{n}H_{\mu^{*}}\left(\left({\displaystyle \bigvee_{k=0}^{n-1}}T^{-k}\alpha_{p}\vee\mathcal{C}\right)^{*}\Big|\mathcal{C}^{*}\right)\\
 &  & =\frac{1}{n}H_{\mu}\left({\displaystyle \bigvee_{k=0}^{n-1}}T^{-k}\alpha_{p}\Big|\mathcal{C}\right) \quad\mbox{by Lemma~\ref{lem:poisson_cond_entropy}} \\\\
 &  & \leq\mu\left(A\right)h\left(A,\mathcal{B}\cap A,\mu_{A},T_{A}\mid\mathcal{C}\cap A\right)
\end{eqnarray*}
by an earlier computation.

Putting things together, we can conclude that \[
h\left(X^{*},\mathcal{B}^{*},\mu^{*},T_{*}\mid\mathcal{C}^{*}\right)\leq\mu\left(A\right)h\left(A,\mathcal{B}\cap A,\mu_{A},T_{A}\mid\mathcal{C}\cap A\right)\]
 which achieves the proof.
\end{proof}

The following corollary is an immediate consequence of Proposition~\ref{prop:relative_entropy}:

\begin{corollary}\label{cor:relative_poisson}
Let $(X,\mathcal{B},\mu,T)$ be a $\mathbf{II}_{\infty}$ system.
\begin{enumerate}
%\item{
%If there exists some factor for which the Poisson and the Parry entropy are equal, then the Poisson entropy of $(X,\mathcal{B},\mu,T)$ is equal to its Parry entropy.}
\item{
If there exists some factor for which the Poisson and the Krengel entropy are equal, then the Poisson entropy of $(X,\mathcal{B},\mu,T)$ is equal to its Krengel entropy.}
\item{If there exists some extension for which the Poisson and the Krengel entropy are equal and finite, then the Poisson entropy of $(X,\mathcal{B},\mu,T)$ is equal to its Krengel entropy.}
\item{If there exists some factor for which the Poisson entropy is zero, then the Poisson entropy of $(X,\mathcal{B},\mu,T)$ is equal to its Parry entropy.}
\item{If there exists some factor for which the Krengel entropy is zero, then the Krengel entropy of $(X,\mathcal{B},\mu,T)$ is equal to its Parry entropy.}
\end{enumerate}
\end{corollary}
\begin{proof}
The first two points are easy consequences of Proposition~\ref{prop:relative_entropy}.

To prove the third point, observe that if $\mathcal{C}$ is a factor on which
the Poisson entropy is zero then, thanks to Proposition~\ref{prop:relative_entropy},
\[
h\left(T^{*}\right)={\displaystyle \sup_{\mathcal{D}\subset\mathcal{B},T_{*}^{-1}\mathcal{D}^{*}\subset\mathcal{D}^{*}}}H_{\mu}\left(\mathcal{D}^{*}\mid T_{*}^{-1}\mathcal{D}^{*}\vee\mathcal{C}^{*}\right)\]
which equals
\[
{\displaystyle \sup_{\mathcal{D}\subset\mathcal{B},T^{-1}\mathcal{D}\subset\mathcal{D}}}H_{\mu}\left(\mathcal{D}\mid T^{-1}\mathcal{D}\vee\mathcal{C}\right)={\displaystyle \sup_{\mathcal{D}\subset\mathcal{B},T^{-1}\mathcal{D}\subset\mathcal{D}}}H_{\mu}\left(\mathcal{D}\vee\mathcal{C}\mid T^{-1}\left(\mathcal{D}\vee\mathcal{C}\right)\right)\]
by the same Proposition. Therefore, $h\left(T^{*}\right)\leq \hpa\left(T\right)$
and the equality follows since Theorem~\ref{thm:poisson_eq_parry} gives the other inequality.

The last point is proven with similar arguments.
\end{proof}

We point out that the assertion (4) of Corollary~\ref{cor:relative_poisson}, which concerns only Krengel and Parry entropy, is implied by~\cite{danilenko_rudolph07}.
%
%Moreover, Danilenko and Rudolph noticed in \cite{danilenko_rudolph07} that the combination of Theorem~\ref{thm:poisson_le_krnegel} and  Proposition~\ref{prop7.1} implies that Krengel and Parry entropy coincide if there exist factors with arbitrarily small Krengel entropy.
%The same reasoning, with the help of Theorem~\ref{thm:poisson_eq_parry} and Proposition~\ref{prop:relative_poisson}, gives:
%\begin{corollary}\label{cor:arbitrarysmallpoissonentropy}
%If $T$ is of type $\mathbf{II}_{\infty}$ and possesses factors with arbitrarily small Poisson entropy, then Poisson entropy and Parry entropy of $T$ are equal.
%\end{corollary}

\section{Additivity and scaling of Poisson entropy}
We now show that just as with Krengel's entropy, the Poisson entropy of a sum of measures is the sum of the Poisson entropies, and scaling a measure by a positive constant scales the Poisson entropy.
\begin{proposition}
\label{prop:p-entropy-additive}
Suppose $\mu$ and $\nu$ are both $T$-invariant $\sigma$-finite measures on $(X,\mathcal{B})$, and
$t,s  >0$. We have
$$h(X^*,\mathcal{B}^*,(t\mu+s\nu)^*,T_*)=s\cdot h(X^*,\mathcal{B}^*,\mu^*,T_*)+t\cdot h(X^*,\mathcal{B}^*,\nu^*,T_*).$$
\end{proposition}
\begin{proof}
Let
$$
(\widehat{X},\widehat{\B},\lambda,\widehat{T}) :=
\Bigl(X\times\{0,1\},\ \mathcal{B}\times 2^{\{0,1\}},\ \mu\times 1_{\{0\}} + \nu \times 1_{\{1\}},\ T \times
Id\Bigr).
$$
This system is isomorphic to the disjoint union of the two
systems $(X,\mathcal{B},\mu,T)$ and $(X,\mathcal{B},\nu,T)$.
The Poisson suspension of $(\widehat{X},\widehat{\B},\lambda,\widehat{T})$ is isomorphic to
the product of the suspensions of $(X,\mathcal{B},\mu,T)$ and $(X,\mathcal{B},\nu,T)$.
Thus,
$$h(\widehat{X}^*,\widehat{\B}^*,\lambda^*,\widehat{T}_*)=h(X^*,\mathcal{B}^*,\mu^*,T_*)+h(X^*,\mathcal{B}^*,\nu^*,T_*).$$
Also, since
$(X,\mathcal{B},\mu+\nu,T)$ is a factor of
$(\widehat{X},\widehat{\B},\lambda,\widehat{T})$, we have that the
suspension of $(X,\mathcal{B},\mu+\nu,T)$ is a factor of the
suspension of $(\widehat{X},\widehat{\B},\lambda,\widehat{T})$. We
thus see that
$$h(X^*,\mathcal{B}^*,(\mu+\nu)^*,T_*)\le h(X^*,\mathcal{B}^*,\mu^*,T_*)+h(X^*,\mathcal{B}^*,\nu^*,T_*).$$
To prove that the above inequality is actually an equality, we
observe that $(\widehat{X},\widehat{\B},\lambda,\widehat{T})$ is a
bounded-to-one extension of $(X,\mathcal{B},\mu+\nu,T)$, and is
therefore a zero-entropy extension.
By Proposition~\ref{prop:relative_entropy}, it follows that
$(\widehat{X}^*,\widehat{\B}^*,\lambda^*,\widehat{T}_*)$ is a zero-entropy extension of $(X^*,\mathcal{B}^*,(\mu+\nu)^*,T_*)$.

We have just proved that Poisson entropy is additive and it remains to
prove that, for any $t>0$,
\begin{equation}\label{eq:poisson_entropy_scaling}
h(X^*,\mathcal{B}^*,(t \cdot \mu)^*,T_*)= t \cdot h(X^*,\mathcal{B}^*, \mu^*,T_*).
\end{equation}
For rational $t$'s, this follows from the above additivity property.
If $t_1 < t_2$, writing $t_2 \cdot\mu = t_1 \cdot \mu + (t_2 - t_1)\cdot \mu$, we get
$h(X^*,\mathcal{B}^*,(t_2 \cdot \mu)^*,T_*)\ge h(X^*,\mathcal{B}^*,(t_1 \cdot \mu)^*,T_*)$.
So  $t \to h(X^*,\mathcal{B}^*,(t \cdot \mu)^*,T_*)$ is a monotone
increasing function. Equation \eqref{eq:poisson_entropy_scaling} now
follows for any real $t>0$, since a monotone function which vanishes
on the rational numbers is zero.
\end{proof}

This result allows to prove that the Poisson entropy of a squashable
transformation is zero or infinite, just as Krengel and Parry entropy (recall that $\left(X,\mathcal{B},\mu,T\right)$ is \emph{squashable} if
it is isomorphic to $\left(X,\mathcal{B},c\mu,T\right)$ for a positive
number $c\neq1$ and \emph{completly squashable} if this holds for
any positive number $c$).

It has been conjectured that stochastic $\alpha$-semi-stable stationary
processes have zero or infinite entropy. It is known in the case $\alpha=2$
which is the Gaussian case (see~\cite{delarue93entgaus}) but remains
unknown otherwise.

However, it has been noticed in~\cite{Roy06IDstat} that $\alpha$-semi-stable stationary
processes ($\alpha<2$) are factors of Poisson suspensions built over
squashable systems (completely squashable in the stable case), associated with the L\'evy measure of the process, which hence are of zero or infinite entropy.

\section{Quasi-finite conservative transformations}
\subsection{Equality of the entropies}
Recall the definition of a quasi-finite transformation from~\cite{krengel67} (also see~\cite{aaro_park}):
Let $(X,\mathcal{B},\mu,T)$ be conservative measure-preserving.
$A \in \mathcal{F}$ is a \emph{quasi-finite set} if $H_{\mu}(\rho_A)< \infty$, where
$\rho_A$ is the first-return-time partition of $A$:
$$
\rho_A:=\left\{A \cap \left(T^{-n}A \setminus\bigcup_{k=1}^{n-1}T^{-k}A\right) ,\ n \ge 1\right\}.
$$
$(X,\mathcal{B},\mu,T)$ is \emph{quasi-finite} if there exists a quasi-finite sweep-out set $A \in \mathcal{F}$.
Using terminology similar to Aaronson and Park~\cite{aaro_park}, we say that a local
partition $\alpha$ is quasi-finite if it has a quasi-finite core $A$, and $\rho_A \prec \alpha$. %This is not in Aaronson and Park.
We point out that conservative transformations which are not quasi-finite have been constructed in~\cite{aaro_park}.

Parry has proved that, for quasi-finite transformations (called ``pseudo-finite'' in~\cite{parry_generators69}), Krengel's definition of entropy coincides with Parry's. We show that, for such transformations, both are equal to the Poisson entropy.

\begin{theorem}\label{thm:poisson_eq_krengel_quasi_finite}
Let $(X,\mathcal{B},\mu,T)$ be a quasi-finite measure-preserving transformation of type $\mathbf{II}_{\infty}$. % and $\alpha$ a quasi
%finite-partition such that
%$$\sup_{n \ge 1} I(\alpha \mid \bigvee_{k=1}^{n}T^{-k}\alpha) \in
%L_1(X,\mathcal{B},\mu)$$
The Poisson entropy, the Krengel entropy and the Parry entropy of
$(X,\mathcal{B},\mu,T)$ are equal. %:
%$$h(X^*,\mathcal{B}^*,\mu^*,T_*) = \hkr (X,\mathcal{B},\mu,T)$$
\end{theorem}

\begin{proof}%[Proof of Theorem~\ref{thm:poisson_eq_krengel_quasi_finite}]
Let $\alpha$ be a local quasi-finite partition whose core $A$ is a sweep-out set.
%Define $\tilde{\alpha}:=\alpha\vee\rho_{A}$.
Applying Proposition~\ref{prop:hat_ge_poisson}, we get
$$
h(X^*, (\hat{\alpha})^*, \mu^*, T_*) \le \hat{h}\left(X,\mu,T,{\alpha}\right).
$$
% and note that even though ${\alpha}$ is not finite, we can replace
% $\alpha$ by ${\alpha}$ in Proposition~\ref{prop:hat_ge_poisson} since $\bigvee_{n\geq0}T^{-n}\alpha=\bigvee_{n\geq0}T^{-n}{\alpha}$,
% $H\left({\alpha}\right)<\infty$ and $\hat{h}\left(X,\mu,T,\alpha\right)\leq\hat{h}\left(X,\mu,T,{\alpha}\right)$.
We want to show that  $\hat{h}\left(X,\mu,T,{\alpha}\right)=H_{\mu}\left(\alpha_{0}^{\infty}\mid\alpha_{1}^{\infty}\right)$.
The result follows by integrating \eqref{information_decomposition} and by proving
the convergence of ${\displaystyle \int_{X}}I_{\mu}\left(\alpha\mid\alpha_{1}^{n}\right)d\mu$ to ${\displaystyle \int_{X}}I_{\mu}\left(\alpha\mid\alpha_{1}^{\infty}\right)d\mu=H_{\mu}\left(\alpha_{0}^{\infty}\mid\alpha_{1}^{\infty}\right)$.
Remark that, since $\rho_{A}\prec{\alpha}$, $I_{\mu}\left(\alpha\mid\alpha_{1}^{n}\right)=I_{\mu}\left(\alpha\mid\alpha_{1}^{\infty}\right)=0$ on $X\setminus A$, therefore
\begin{equation}
\label{eq:return_time}
\int_{X}I_{\mu}\left(\alpha\mid\alpha_{1}^{n}\right)d\mu = \int_{A} I_{\mu}\left(\alpha\mid\alpha_{1}^{n}\right)d\mu.
\end{equation}
By Lemma~\ref{lem:info_mean_converge} applied to the set $A$ with the restriction of the
$\sigma$-algebras ${\alpha}$ and $\alpha_{1}^{n}$ to $A$, the right-hand side tends to ${\displaystyle \int_{A}}I_{\mu}\left(\alpha\mid\alpha_{1}^{\infty}\right)d\mu={\displaystyle \int_{X}}I_{\mu}\left(\alpha\mid\alpha_{1}^{\infty}\right)d\mu$
which gives us the desired convergence.

Putting things together, we have proved
\[
h\left(X^{*},\left(\alpha_{0}^{\infty}\right)^{*},\mu^{*},T_{*}\right)\le H_{\mu}\left(\alpha_{0}^{\infty}\mid\alpha_{1}^{\infty}\right),\]
the right hand-side being bounded by $\hpa\left(X,\widehat{{\alpha}},\mu,T\right)$
by definition.
By Theorem~\ref{thm:poisson_eq_parry}, the latter is in turn dominated by $h\left(X^{*},\left(\widehat{{\alpha}}\right)^{*},\mu^{*},T_{*}\right)$. Hence, we obtain
\begin{equation}
h\left(X^{*},\left(\widehat{{\alpha}}\right)^{*},\mu^{*},T_{*}\right)=H_{\mu}\left(\alpha_{0}^{\infty}\mid\alpha_{1}^{\infty}\right)=\hpa\left(X,\widehat{{\alpha}},\mu,T\right).
\label{eq:EDpartition}
\end{equation}
Now replace $\alpha$ by $\alpha_{n}$ in \eqref{eq:EDpartition},
where $(\alpha_{n})$ is an increasing sequence of local quasi-finite partitions
with core $A$ having the property that $\left(\alpha_{n}\right)_{0}^{\infty}\uparrow\mathcal{B}$.
Taking the limit in $n$, we obtain
$$
h\left(X^{*},\mathcal{B}^{*},\mu^{*},T_{*}\right)=\hpa\left(X,\mathcal{B},\mu,T\right),
$$
\textit{i.e.} the Poisson entropy equals the Parry entropy. At last, the Krengel entropy equals
the two others since the system is quasi-finite.
\end{proof}

% We define a measurable partition $\xi$ to be \emph{entropy determining} (abbreviated ED) if it contains at least one non-empty set of finite measure,
% $H_{\mu}(\xi) < \infty$, and $H_{\mu}(\xi_{\infty}^0 \mid T^{-1}\xi_{\infty}^0)$ is equal to the Parry entropy of $T$ with respect to the $\sigma$-algebra $\hat\xi$.
%
% In the course of proof of Theorem~\ref{thm:poisson_eq_krengel_quasi_finite}, we have actually established the following:
% \begin{proposition}\label{prop:qf_is_EF}
% Any quasi-finite partition is ED.
% \end{proposition}

\subsection{Poisson suspensions of Markov chains}
Poisson suspensions of Markov chains have been considered by several authors.
Grabinsky~\cite{grabinsky84} and Kalikow~\cite{kalikow81} have independently proved that the Poisson suspension of an ergodic, null-recurrent random walk is Bernoulli.

Let $\Sigma$ be a countable or finite set, $P= (p_{a,b})_{a,b \in
\Sigma}$ be a stochastic matrix which is irreducible and recurrent.
As is well-known, there exists a measure $q$ on $\Sigma$
which is stationary with respect to $P$, meaning $qP=q$, and this
measure is unique up to scaling.
The associated Markov shift is the system $(X,\mathcal{B},\mu,T)$, where $X=\Sigma^\mathbb{Z}$, $T:X
\to X$ denotes the shift map ($T(x)_n=x_{n+1}$), $\mathcal{B}$ denotes
the Borel $\sigma$-algebra of $X$ with respect to the product
topology and $\mu$ is given by
$$\mu([a_1,\ldots,a_k])=q_{a_1}\prod_{i=2}^k p_{a_{i-1},a_i}.$$

Based on the Krengel entropy of recurrent Markov chains and our previous result about Poisson entropy of quasi-finite transformations, we have
\begin{corollary}
\label{cor:markov_entropy}
The entropy of the Poisson suspension of a recurrent Markov shift with transition matrix $P= (p_{a,b})_{a,b \in
\Sigma}$ and stationary measure $q$ is given by
\begin{equation}\label{eq:markov_entropy}
h(X^*,\mathcal{B}^*,\mu^*,T^*)=\sum_{a \in \Sigma}q_a \sum_{b \in \Sigma}p_{a,b}\log\frac{1}{p_{a,b}}.
\end{equation}
\end{corollary}

\begin{proof}
By Krengel's formula (Theorem~$4.1$ of~\cite{krengel67}), the Krengel entropy of $(X,\mathcal{B},\mu,T)$ is given by the right-hand side of~\eqref{eq:markov_entropy}. By taking the standard Markov partition $\xi$, we see that
$$H(\xi^0_{-\infty} \mid T^{-1}\xi^0_{-\infty})= \sum_{a \in \Sigma}q_a \sum_{b \in \Sigma}p_{a,b}\log\frac{1}{p_{a,b}}.$$
Thus, Parry's entropy dominates Krengel's. Hence both are equal.

Without loss of generality we can assume that the transition matrix is irreducible.
In the particular case when $(X,\mathcal{B},\mu,T)$ is a renewal system ($\Sigma=\NN$ and $p_{n,n-1}=1$ for all $n >1$), the right-hand side of \eqref{eq:markov_entropy} is simply the entropy of the first-return-time partition of the state $1$. Hence, if the Krengel entropy is finite, the renewal system is quasi-finite, in which case the Poisson entropy is equal to the Krengel entropy by Theorem~\ref{thm:poisson_eq_krengel_quasi_finite}. Otherwise, the Parry entropy is infinite, and by Theorem~\ref{thm:poisson_eq_parry} it is equal to the Poisson entropy.
Now it remains to note that every irreducible recurrent Markov chain has a factor which is a renewal system.
Hence, a Markov chain is quasi-finite if and only if it has finite Krengel entropy, in which case this is also the Poisson entropy. Otherwise, the Poisson entropy is infinite.
\end{proof}

In~\cite{grabinsky84}, it is claimed that the entropy of the Poisson suspension of a null-recurrent Markov chain is infinite (Proposition 4.3).
Corollary~\ref{cor:markov_entropy} together with the existence of such chains with finite Krengel entropy (see~\cite{krengel67}) contradict this result.
The mistake in~\cite{grabinsky84} comes from the following incorrect assertion which Grabinsky invokes in the proof of Proposition 4.3, to bound from below the entropy of a certain partition: If $\xi$ is a Markov partition (i.e. $\xi$ is independent of $\xi^{-1}_{-\infty}$ given $T^{-1}\xi$) and $\eta$ a partition which is measurable with respect to $\xi$, then $\eta^*$ is a Markov partition as well. It would imply, in particular, that given the number of
particles in a certain Markov state $A$ at time $-1$, the number of particles in $A$ at  time 0
is independent of the number of particles in $A$ at time $-2$\dots

%I found the following mistakes in Grabinsky's proof:
%on page 310, line -3, it is claimed that if Q is a partition which is measurable with
%respect to the "Poissonification" of the last r-states,
%then given the recent r-states of Q, it is independent of the previous history. This
%assertion is in correct: It would imply, in particular, that given the number of
%particles in a certain Markov state A at time -1, the number of particles in A at  time 0
%is independent of the number of particles in A at time -2...
%Grabinsky claims that this follows from his lemma 3.1 - but it doesn't, as a
%sub-partition of a Markov partition is not always Markov.
%
%Then in the proof of proposition 4.3 (the one which we contradict)
%P. 312, line 7
%it is claimed that
%H_p(T_u) \ge \sup\{ H_p(Q | T_u^{-1} Q): \sigma(Q) \subset \gamma_0, H_p(Q) < \infty\}
%
%This is false, for the same reason: H_p(Q | T_u^{-1} Q) should be replaced by H_p(Q |
%\bigvee_{n > 0} T_u^{-n}Q), which  can't be bounded from below by something that tends to
%\infty...

\section{\label{sec:zero_entropy} Zero Poisson entropy}

In this section, we prove that some class of cutting-and-stacking constructions (including finite-rank transformations) and transformations without a countable Lebesgue  component in their spectrum both have zero Poisson entropy. It is well known that  rank one transformations also have zero Krengel entropy, therefore, Krengel entropy, Parry entropy and Poisson entropy are equal in this case.

The construction of a non-quasi-finite transformation in~\cite{aaro_park} is of this kind, so the results of this section do not follow from Theorem~\ref{thm:poisson_eq_krengel_quasi_finite}.

\subsection{Cutting-and-stacking constructions}
A cutting-and-stacking construction is an iterative method to present conservative transformations.
We briefly describe this construction and refer to Friedman's book~\cite{friedman_book} for details.

A \emph{column} of height $h \in \mathbb{N}$ is an array $I_1,\ldots,I_h$ of pairwise disjoint intervals of the same length, considered as ``stacked'' one on top of the other.
At stage $n$ of the cutting-and-stacking procedure, the $n$-th \emph{tower} $X_n$ consists of $c_n$ columns of heights $\{h_{n,i}\}_{1\le i \le c_n}$ and equal width.
The transformation acts by translating each interval to the interval above it.
At stage $n$, the transformation is undefined for points on the top intervals.
At stage $n+1$, each column is ``cut'' into $k_n$ columns, all the columns are ``stacked'' one on top of the other, and then the newly formed column is cut into $c_{n+1}$ columns.
Then some new intervals are possibly added on the top of every column.
As the length of the intervals at stage $n$ tends to $0$, the measure of the points on which the transformation is undefined at stage $n$ tends to $0$.
Such a construction is said to have \emph{rank one} if $c_n=1$ for every $n\ge 1$, and \emph{finite rank} if $\{c_n\}$ is bounded.

Denote by $\epsilon_n$ the length of the intervals at stage $n$. Clearly, $\epsilon_n \to 0$ as $n \to \infty$.

\begin{proposition}
\label{prop:cutting_and_stacking_zero_entropy}
Let $(X,\mathcal{B},\mu,T)$ be a cutting-and-stacking construction as above.
If $c_n \epsilon_n \log \epsilon_n \to 0$ as $n\to \infty$, then $h(X^*,\mathcal{B}^*,\mu^*,T_*)=0$.
In particular, this is the case if $T$ has finite rank.
\end{proposition}

\begin{proof}
Let $\beta_n= \{I_{n,1},\ldots,I_{n,c_n}\}$ denote the set of intervals composing the base of the tower at stage $n$, and let
$$
\xi_n =\{I_{n,1}, T I_{n,1},\ldots,T^{h_{n,1}-1}I_{n,1},\ldots,I_{n,c_n},\ldots,T^{h_{n,c_n}-1}I_{n,c_n}\}
$$
denote the corresponding partition of $X_n$.
We have $\xi_n^* \subset \hat{\beta_n^*}$. Since $\xi_n \uparrow \mathcal{B}$, $\hat{\beta_n^*} \uparrow \mathcal{B}^*$.
{F}rom this, we deduce that $h_{\mu^*}(T_*,\hat{\beta_n^*}) \uparrow h(X^*,\mathcal{B}^*,\mu^*,T_*)$.
But $h_{\mu^*}(T_*,\hat{\beta_n^*}) \le H_{\mu^*}(\beta_n^*) = c_n f(\epsilon_n) \sim -c_n \epsilon_n \log \epsilon_n$, where $f(x)$ denotes the entropy function of a Poisson variable with parameter $x$.
\end{proof}

%A certain generalization of rank-one cutting-and-stacking constructions into countable group actions is known as $(C,F)$-contractions \cite{danilenko_rank_one}. The above proof of Proposition~\ref{prop:cutting_and_stacking_zero_entropy} actually gives:
%
%\begin{proposition}
%A $(C,F)$-construction over any countable amenable group has zero Poisson entropy.
%\end{proposition}

\subsection{Spectral criterion}

The following proposition and corollary give a spectral criterion for positive Poisson entropy. The corresponding result about Kolmogorov entropy is well known in the finite measure case.

\begin{proposition}\label{pro:countable_lebesgue}
If $\left(X,\mathcal{B},\mu,T\right)$ is of type $\mathbf{II}_{\infty}$
and has positive Poisson entropy, then its spectrum has a component
which is countable Lebesgue.
\end{proposition}
\begin{proof}
Pick a sweep out set $A\in\mathcal{B}$ with small measure in order
to have $H_{\mu^*}\left(\xi_{0}^{*}\right)=f\left(\mu\left(A\right)\right)<h\left(T_{*}\right)$,
where $\xi_{0}$ is the partition $\left\{ A,A^{c}\right\} $ of $X$.
We now refine the local partition $\xi_{0}$ by increasing finite local
partitions $\xi_{n}$ so that $\xi_{n}\uparrow\mathcal{B}$ which
implies $h\left(T_{*},\xi_{n}^{*}\right)\uparrow h\left(T_{*}\right)$.
Using the continuity of Kolmogorov entropy, the continuity of $f$
and the continuity of $\mu$, we can assume that this increasing sequence
$\xi_{n}$ is such that
\[
0<h\left(T_{*},\xi_{1}^{*}\right)<\cdots<h\left(T_{*},\xi_{n}^{*}\right)<\cdots\]

In the following, for a $\sigma$-algebra $\mathcal{A}\subset\mathcal{B}^{*}$,
we denote by $L^{2}\left(\mathcal{A}\right)$ the corresponding linear
subspace of $L^{2}\left(\mu^{*}\right)$ of square-integrable $\mathcal{A}$-measurable
functions, and $U_{T_{*}}$ is the unitary operator induced from $T_{*}$. Denote by $\mathfrak{C}$, the \emph{first chaos} of $L^{2}\left(\mu^{*}\right)$,
i.e. the closure of the linear subspace generated by $N\left(A\right)-\mu\left(A\right)$,
$A\in\mathcal{B}$, $\mu\left(A\right)<\infty$. The arguments below are classical; we already know that the suspension has a countable Lebesgue component in its spectrum, however, to get the result, we will see that a countable Lebesgue component is localized in $\mathfrak{C}$ which is unitary isomorphic to $L^{2}\left(\mu\right)$.
 Set $H_{1}:=L^{2}\left(\left(\xi_{1}^{*}\right)_{-\infty}^{0}\right)\cap\mathfrak{C}$,
and note that it is non-empty since it contains the functions $N\left(A\right)-\mu\left(A\right),A\in T^{-k}\xi_{1}$,
$k\in\mathbb{N}$. Remark that $U_{T_{*}^{-1}}H_{1}\subset H_{1}$
and that we cannot have $U_{T_{*}^{-1}}H_{1}=H_{1}$ since it would
imply that $\sigma\left(H_{1}\right)$ belongs to the Pinsker factor
of $T_{*}$, and, as the factor $\widehat{\xi_{1}^{*}}$ is generated
(as $\sigma$-algebra) by $\cup_{n\geq0}U_{T_{*}^{n}}H_{1}=H_{1}$,
$h\left(T_{*},\xi_{1}^{*}\right)=0$ which is a contradiction. Functions
belonging to $V_{1}:=H_{1}\ominus U_{T_{*}^{-1}}H_{1}$ have Lebesgue
spectral measure.

Set $H_{2}:=\left(L^{2}\left(\left(\xi_{2}^{*}\right)_{-\infty}^{0}\right)\cap\mathfrak{C}\right)\cap\left(L^{2}\left(\widehat{\xi_{1}^{*}}\right)\cap\mathfrak{C}\right)^{\perp}$.
It is also non-empty since $\cup_{n\geq0}U_{T_{*}^{n}}\left(L^{2}\left(\left(\xi_{2}^{*}\right)_{-\infty}^{0}\right)\cap\mathfrak{C}\right)$
generates $\widehat{\xi_{2}^{*}}$, which is strictly bigger than $\widehat{\xi_{1}^{*}}$
because $h\left(T_{*},\xi_{1}^{*}\right)<h\left(T_{*},\xi_{2}^{*}\right)$.
Moreover, we have $U_{T_{*}^{-1}}H_{2}\subset H_{2}$, and once again,
we cannot have $U_{T_{*}^{-1}}H_{2}=H_{2}$ as it would imply that
$\sigma\left(H_{2}\right)$ belongs to the Pinsker factor of $T_{*}$,
and this is impossible since the entropy of $\widehat{\xi_{2}^{*}}=\widehat{\xi_{1}^{*}}\vee\sigma\left(H_{2}\right)$
would be equal to that of $\widehat{\xi_{1}^{*}}$ which is a contradiction.
Therefore we can set $V_{2}:=H_{2}\ominus U_{T_{*}^{-1}}H_{2}$, which
is constituted by functions having Lebesgue spectral measure and satisfies
$\overline{\cup_{n\in\mathbb{Z}}U_{T_{*}^{n}}V_{2}}\perp\overline{\cup_{n\in\mathbb{Z}}U_{T_{*}^{n}}V_{1}}$.

Proceeding by induction, we construct an infinite sequence of mutually
orthogonal invariant subspaces $\overline{\cup_{n\in\mathbb{Z}}U_{T_{*}^{n}}V_{n}}$ of $\mathfrak{C}$
on which $U_{T_{*}}$ has Lebesgue maximal spectral type. Thanks to
the unitary isomorphism between $\mathfrak{C}$ and $L^{2}\left(\mu\right)$,
these subspaces can be transferred into $L^{2}\left(\mu\right)$ and
we have proved that $T$ has a countable Lebesgue component in its
spectrum.
\end{proof}

\begin{corollary}\label{cor:singular_spectrum}
If the maximal spectral type of $T$ is singular or if T has finite multiplicity, then its Poisson entropy is zero.
\end{corollary}

We mention that Parry~\cite{parry_spectral} has shown that a $K$-automorphism has countable Lebesgue spectrum.
{F}rom his proof one can easily obtain that if the maximal spectral type is singular or has finite multiplicity, then the Parry entropy is zero.
Since zero Poisson entropy implies zero Parry entropy, our above corollary refines Parry's result.
We do not know of a sufficient spectral criterion for zero Krengel entropy.

\section{\label{sec:perfect_partitions} Perfect Poissonian partitions and the Poisson-Pinsker factor}
For a probability-preserving system $(X,\mathcal{B},\mu,T)$, the \emph{Pinsker factor}, denoted by $\mathcal{P}(T)$, is the maximum factor ($T$-sub-invariant $\sigma$-algebra) with zero entropy.
%By a theorem of Rokhlin and Sina\"{i}, this factor is made of sets constituting tail $\sigma$-algebra of all partitions of finite entropy.
For each of the various notions of entropy for $\sigma$-finite transformations, we generalize this definition: We say that a factor is Pinsker if it is the maximum zero entropy factor.
We can thus speak of a Krengel-Pinsker factor, a Parry-Pinsker factor and a Poisson-Pinsker factor of a conservative transformation.
The existence of these is not obvious in general.
The following proposition gives a sufficient condition for Pinsker factors to exist, which we later show to be necessary as well:

\begin{proposition}\label{prop:PoissonPinsker}
Let $(X,\mathcal{B},\mu,T)$ be an ergodic type $\mathbf{II}_{\infty}$ system. Assume the Pinsker factor $\mathcal{P}\left(T_{*}\right)$
of the Poisson suspension is of the form $\mathcal{P}^{*}$ for some $\sigma$-finite
$\sigma$-algebra $\mathcal{P}$. Then $\mathcal{P}$ is both the Poisson-Pinsker
and the Parry-Pinsker factor of $T$. Moreover, if there exists a factor with zero Krengel
entropy, then $\mathcal{P}$ is also the Krengel-Pinsker factor.
\end{proposition}
\begin{proof}
Since $\mathcal{P}^*$ is the Pinsker  factor of $T_*$, for any factor $\mathcal{C}$ of zero Poisson entropy, $\mathcal{C}^* \subset \mathcal{P}^*$. Since $\mathcal{C}$ is $\sigma$-finite, this implies that $\mathcal{C} \subset \mathcal{P}$. This shows that $\mathcal{P}$ is the Poisson-Pinsker factor of $(X,\mathcal{B},\mu,T)$.

We now prove that $\mathcal{P}$ is also the Parry-Pinsker factor.
Assume that $\mathcal{C}$ is a factor of zero Parry entropy.
By Lemma~\ref{lem:star-intersection}, we have $\mathcal{P}^{*}\cap\mathcal{C}^{*}=\left(\mathcal{P}\cap\mathcal{C}\right)^{*}$.
Moreover, if $\mathcal{P}\cap\mathcal{C}$ contains no sets of positive finite measure,
then $\left(\mathcal{P}\cap\mathcal{C}\right)^{*}$ is trivial, and since $\mathcal{P}^{*}$
is the Pinsker factor of $\left(X^{*},\mathcal{B}^{*},\mu^{*},T_{*}\right)$,
then $\mathcal{C}^{*}$ is a $K$-system.
By a well-known disjointness result (for probability-preserving transformations),
$\mathcal{P}^{*}$ and $\mathcal{C}^{*}$ are independent. For any
two positive and finite measure sets $A\in\mathcal{P}$ and $B\in\mathcal{C}$,
we have
\[
{\displaystyle \int_{X^{*}}}\Bigl(\gamma\left(A\right)-\mu\left(A\right)\Bigr)\Bigl(\gamma\left(B\right)-\mu\left(B\right)\Bigr)\mu^{*}\left(d\gamma\right)=0.\]
But the left-hand side equals $\mu\left(A\cap B\right)$, so $A$ and $B$ are disjoint mod. $\mu$.
This is impossible for all $A\in\mathcal{P}$ and $B\in\mathcal{C}$ because it
would contradict the fact that $\mathcal{P}$ is $\sigma$-finite
and $T$ of type $\mathbf{II}_{\infty}$. This shows that $\mathcal{P}\cap\mathcal{C}$ must contain a set of positive, finite measure, so by ergodicity it is $\sigma$-finite. Thus, $\mathcal{P}\cap\mathcal{C}$ is a factor of zero Poisson entropy of $\mathcal{C}$. Thanks to Corollary~\ref{cor:relative_poisson}~(3), Parry and Poisson entropy coincide on $\mathcal{C}$ and since the first one is zero, the second one is zero. This implies that $\mathcal{C}\subset\mathcal{P}$.
Hence $\mathcal{P}$ is also the Parry-Pinsker factor. The statement about
Krengel-Pinsker factor is proved in the same way, using Corollary~\ref{cor:relative_poisson}~(4).
\end{proof}

Recall that a $\sigma$-algebra $\xi$ for an invertible probability-preserving
system $\left(\Omega,\mathcal{F},m,T\right)$ is called a \emph{perfect
$\sigma$-algebra } if $T^{-1}\xi\subset\xi$, $\widehat{\xi}=\mathcal{B}$,
$h\left(T\right)=H_{m}\left(\xi\mid T^{-1}\xi\right)$ and \[
{\displaystyle \bigcap_{n=0}^{\infty}}T^{-n}\xi=\mathcal{P}\left(\Omega,\mathcal{F},m,T\right).\]

In fact, if the entropy is finite, the last condition is a consequence of the others as it is
proved in the following lemma.
The ingredients are similar to the proof of Rokhlin-Sina\"{i}  Theorem,
as appearing in~\cite{parry_topics}, page 69.

\begin{lemma}\label{lem:perfect_partitions}
Assume $\left(\Omega,\mathcal{F},m,T\right)$ is an invertible probability-preserving
system. If $T^{-1}\xi\subset\xi$, $\widehat{\xi}=\mathcal{B}$ and $h\left(T\right)=H_{m}\left(\xi\mid T^{-1}\xi\right)<\infty$,
then \[
{\displaystyle \bigcap_{n=0}^{\infty}}T^{-n}\xi=\mathcal{P}\left(\Omega,\mathcal{F},m,T\right).\]
\end{lemma}
\begin{proof}
Let $\alpha_{n}$ be a sequence of finite entropy partitions increasing
to $\xi$. Since $h\left(T\right)=H_{m}\left(\xi\mid T^{-1}\xi\right)$,
we have
\[
h\left(T,\widehat{\alpha_{n}}\right)=H_{m}\left(\alpha_{n}\mid\left(\alpha_{n}\right)_{-\infty}^{-1}\right)\uparrow h\left(T\right).
\]
Since $H_{m}\left(\alpha_{n}\mid T^{-1}\xi\right)$ also converges to $h(T)$, passing to a subsequence, we can assume
\[
H_{m}\left(\alpha_{n}\mid\left(\alpha_{n}\right)_{-\infty}^{-1}\right)-H_{m}\left(\alpha_{n}\mid T^{-1}\xi\right)\leq\frac{1}{n}.
\]
Let $\zeta\subset{\displaystyle \bigcap_{n=0}^{\infty}}T^{-n}\xi$.
It follows that $\widehat{\zeta}\subset{\displaystyle \bigcap_{n=0}^{\infty}}T^{-n}\xi$.
Applying the formula (Theorem 8, page 66 in~\cite{parry_topics}),
\[
h\left(T,\beta\vee\zeta\right)=h\left(T,\zeta\right)+H_{m}\left(\beta\mid\widehat{\zeta}\vee\beta_{-\infty}^{-1}\right)\]
we obtain
\[
H_{m}\left(\zeta\mid\zeta_{-\infty}^{-1}\right)=H_{m}\left(\zeta\vee\alpha_{n}\mid\zeta_{-\infty}^{-1}\vee\left(\alpha_{n}\right)_{-\infty}^{-1}\right)-H_{m}\left(\alpha_{n}\mid\widehat{\zeta}\vee\left(\alpha_{n}\right)_{-\infty}^{-1}\right)\]
$\leq H_{m}\left(\alpha_{n}\mid\left(\alpha_{n}\right)_{-\infty}^{-1}\right)+H_{m}\left(\zeta\mid\left(\alpha_{n}\right)_{-\infty}^{-1}\right)-H_{m}\left(\alpha_{n}\mid T^{-1}\xi\right)\leq\frac{1}{n}+H_{m}\left(\zeta\mid\left(\alpha_{n}\right)_{-\infty}^{-1}\right)$,
which goes to zero as $n$ tends to infinity.
Therefore, $h\left(T,\widehat{\zeta}\right)=0$
and $\zeta\in\mathcal{P}\left(\Omega,\mathcal{F},m,T\right)$.
We thus proved that ${\displaystyle \bigcap_{n=0}^{\infty}}T^{-n}\xi\subset\mathcal{P}\left(\Omega,\mathcal{F},m,T\right)$.
The other inclusion is a consequence of Theorem 13 page 69 in~\cite{parry_topics}, which states
that for any increasing sequence of strictly invariant $\sigma$-algebras
$\mathcal{B}_{n}\uparrow\mathcal{F}$, then
\[
\mathcal{B}_{n}\cap\mathcal{P}\left(\Omega,\mathcal{F},m,T\right)\uparrow\mathcal{P}\left(\Omega,\mathcal{F},m,T\right).\]
Indeed, applying this result with $\mathcal{B}_{n}:=\widehat{\alpha_{n}}$, we get that \[
\widehat{\alpha_{n}}\cap\mathcal{P}\left(\Omega,\mathcal{F},m,T\right)\uparrow\mathcal{P}\left(\Omega,\mathcal{F},m,T\right).\]
Moreover, since $\alpha_{n}$ is a finite entropy partition, \[
\widehat{\alpha_{n}}\cap\mathcal{P}\left(\Omega,\mathcal{F},m,T\right)={\displaystyle \bigcap_{k=0}^{\infty}}\left(\alpha_{n}\right)_{-\infty}^{-k},\]
which is included in ${\displaystyle \bigcap_{j=0}^{\infty}}T^{-j}\xi$
because $\alpha_{n}\subset\xi$ for any $n$.
It follows that $\mathcal{P}\left(\Omega,\mathcal{F},m,T\right)\subset{\displaystyle \bigcap_{n=0}^{\infty}}T^{-n}\xi$.
\end{proof}

Let us introduce the following definition:
A $\sigma$-finite $\sigma$-algebra $\mathcal{A}$ is said to be \emph{entropy determining (ED)} if $T^{-1}\mathcal{A}\subset \mathcal{A}$ and $\mathcal{A}^*$ is perfect with respect to the factor it generates.
Observe that on the factor generated by an ED $\sigma$-algebra $\mathcal{A}$, Parry and Poisson entropies coincide:
$$\hpa(T,\widehat{\mathcal{A}}) \le h(T_*,\mathcal{A}^*) = H_{\mu^*}(\mathcal{A}^* | T_*^{-1}\mathcal{A}^* )
=H_{\mu}(\mathcal{A} | T^{-1}\mathcal{A}) \le \hpa(T,\widehat{\mathcal{A}}).$$
The class of ED $\sigma$-algebras plays the same role as finite-entropy partitions do in the finite-measure case.

We are now ready to prove a ``Poisson analogue'' of the Rokhlin-Sina\"{i}  Theorem regarding Pinsker factors of probability-preserving transformations.
We recall that $T$ is \emph{remotely infinite} if there exists a $\sigma$-finite sub-$\sigma$-algebra $\alpha$ such that $T^{-1}\alpha\subset\alpha$, $T^{n}\alpha\uparrow\mathcal{B}$ and $T^{-n}\alpha\downarrow\mathcal{G}$ (mod. $\mu$) where $\mathcal{G}$ has no set of positive finite measure.
%If $T$ is remotely infinite and $\mathcal{G}=\left\{ X,\emptyset\right\}$, $T$ is called a \emph{$K$-system}.
Also recall that a probability-preserving transformation is a $K$-system if and only if there exists a sub-invariant generating $\sigma$-algebra with a trivial tail.
The notion of remotely-infinite system can be considered as an infinite-measure analogue of a probability-preserving $K$-system. The Rokhlin-Sina\"{i} Theorem tells us that a probability-preserving transformation is a $K$-system if and only if the trivial factor is the only factor of zero entropy.
\begin{theorem}\label{thm:poisson_rohlin_sinai}
Let $(X,\mathcal{B},\mu,T)$ be an ergodic system of type $\mathbf{II}_{\infty}$. Assume there exists an ED partition $\mathcal{A}$ such that $H_{\mu}(\mathcal{A}\mid T^{-1}\mathcal{A})<\infty$. Then
\begin{itemize}
\item there exists a generating ED partition (in particular $\hpa \left(T\right)=h\left(T_{*}\right)$).
\item $T$ is either  remotely infinite or $\mathcal{P}\left(T_{*}\right)$ is Poissonian:  $\mathcal{P}\left(T_{*}\right)=\mathcal{P}^*$  for some $\sigma$-finite $T$-invariant $\sigma$-algebra $\mathcal{P}$. In the latter case, $\mathcal{P}$ is the Poisson (and Parry) Pinsker factor of $T$.
\end{itemize}
\end{theorem}

\begin{proof}
Let $\xi$ be a finite local partition with a sweep-out core $A\in\mathcal{A}$.
We first show that $\xi_{-\infty}^0\vee\mathcal{A}$ is also an ED $\sigma$-algebra.

On the one hand, we have
\begin{eqnarray*}
 \lefteqn{H_{\mu^{*}}\left(\left(\left(\xi_{-p}^{0}\right)^{*}\vee\mathcal{A}^{*}\right)_{0}^{n}\mid\mathcal{A}^{*}\right)}\\
 & = & H_{\mu^{*}}\left(\left(\left(\xi_{-p}^{0}\right)^{*}\right)_{0}^{n}\mid\left(\mathcal{A}^{*}\right)_{0}^{n}\right)+H_{\mu^{*}}\left(\left(\mathcal{A}^{*}\right)_{0}^{n}\mid\mathcal{A}^{*}\right)\\
 & =  & H_{\mu^{*}}\left(\left(\left(\xi_{-p}^{0}\right)^{*}\right)\mid\mathcal{A}^{*}\right)+{\displaystyle \sum_{k=1}^{n}}H_{\mu^{*}}\left(\left(\xi_{-p}^{0}\right)^{*}\vee\mathcal{A}^{*}\mid T_{*}^{-1}\left(\left(\xi_{-p}^{0}\right)^{*}\right)_{-k}^{0}\vee T_{*}^{k}\mathcal{A}^{*}\right)
\\
&& \hspace{3cm} +\ nH_{\mu^{*}}\left(\mathcal{A}^{*}\mid T_{*}^{-1}\mathcal{A}^{*}\right).
\end{eqnarray*}
Dividing by $n$, this tends to
\[
H_{\mu^{*}}\left(\left(\xi_{-p}^{0}\right)^{*}\mid T_{*}^{-1}\left(\left(\xi_{-p}^{0}\right)^{*}\right)_{-\infty}^{0}\vee\widehat{\mathcal{A}^{*}}\right)+H_{\mu^{*}}\left(\mathcal{A}^{*}\mid T_{*}^{-1}\mathcal{A}^{*}\right).\]
Since $\mathcal{A}$ is ED, the second term equals $h\left(T_{*},\mathcal{A}^{*}\right)$. Thus the previous expression equals
$h\left(T_{*},\left(\xi_{-p}^{0}\right)^{*}\vee\mathcal{A}^{*}\right)$.

On the other hand,
\begin{eqnarray*}
 &  & H_{\mu}\left(\left(\xi_{-p}^{0}\vee\mathcal{A}\right)_{0}^{n}\mid\mathcal{A}\right)\\
 &  & =H_{\mu}\left(\xi_{-p}^{0}\vee\mathcal{A}\mid\mathcal{A}\right)+{\displaystyle \sum_{k=1}^{n}}H_{\mu}\left(T^{k}\left(\xi_{-p}^{0}\vee\mathcal{A}\right)\Bigg\vert\left({\displaystyle \bigvee_{j=0}^{k-1}}T^{j}\left(\xi_{-p}^{0}\vee\mathcal{A}\right)\right)\vee\mathcal{A}\right)\\
 &  & =H_{\mu}\left(\xi_{-p}^{0}\vee\mathcal{A}\mid\mathcal{A}\right)+{\displaystyle \sum_{k=1}^{n}}H_{\mu}\left(T^{k}\left(\xi_{-p}^{0}\vee\mathcal{A}\right)\Bigg\vert\left({\displaystyle \bigvee_{j=0}^{k-1}}T^{j}\xi_{-p}^{0}\right)\vee T^{k-1}\mathcal{A}\right)\\
 &  & =H_{\mu}\left(\xi_{-p}^{0}\vee\mathcal{A}\mid\mathcal{A}\right)+{\displaystyle \sum_{k=1}^{n}}H_{\mu}\left(T^{k}\left(\xi_{-p}^{0}\vee\mathcal{A}\right)\Bigg\vert T^{k}\left(\left({\displaystyle \bigvee_{j=-k}^{-1}}T^{j}\xi_{-p}^{0}\right)\vee T^{-1}\mathcal{A}\right)\right)\\
 &  & =H_{\mu}\left(\xi_{-p}^{0}\vee\mathcal{A}\mid\mathcal{A}\right)+{\displaystyle \sum_{k=1}^{n}}H_{\mu}\left(\xi_{-p}^{0}\vee\mathcal{A}\ \Big\vert\ T^{-1}\left(\xi_{-p}^{0}\right)_{-k}^{0}\vee T^{-1}\mathcal{A}\right),
\end{eqnarray*}
which, divided by $n$, tends to
$H_{\mu}\left(\xi_{-\infty}^{0}\vee\mathcal{A}\mid T^{-1}\left(\xi_{-\infty}^{0}\vee\mathcal{A}\right)\right)$.

But observe that
\begin{eqnarray*}
H_{\mu^{*}}\left(\left(\left(\xi_{-p}^{0}\right)^{*}\vee\mathcal{A}^{*}\right)_{0}^{n}\mid\mathcal{A}^{*}\right) & \leq & H_{\mu^{*}}\left(\left(\left(\xi_{-p}^{0}\vee\mathcal{A}\right)^{*}\right)_{0}^{n}\mid\mathcal{A}^{*}\right)\\
\leq H_{\mu^{*}}\left(\left(\left(\xi_{-p}^{0}\vee\mathcal{A}\right)_{0}^{n}\right)^{*}\mid\mathcal{A}^{*}\right) & = & H_{\mu}\left(\left(\xi_{-p}^{0}\vee\mathcal{A}\right)_{0}^{n}\mid\mathcal{A}\right).
\end{eqnarray*}
Therefore, we have for all $p$
\[
h\left(T_{*},\left(\xi_{-p}^{0}\right)^{*}\vee\mathcal{A}^{*}\right)\leq H_{\mu}\left(\xi_{-\infty}^{0}\vee\mathcal{A}\mid T^{-1}\left(\xi_{-\infty}^{0}\vee\mathcal{A}\right)\right)\]
\[
=H_{\mu^{*}}\left(\left(\xi_{-\infty}^{0}\vee\mathcal{A}\right)^{*}\mid T_{*}^{-1}\left(\xi_{-\infty}^{0}\vee\mathcal{A}\right)^{*}\right)\]
and we deduce that
\[
h\left(T_{*},\left(\xi_{-\infty}^{0}\right)^{*}\vee\mathcal{A}^{*}\right)\leq H_{\mu^{*}}\left(\left(\xi_{-\infty}^{0}\vee\mathcal{A}\right)^{*}\mid T_{*}^{-1}\left(\xi_{-\infty}^{0}\vee\mathcal{A}\right)^{*}\right).\]
But $\xi_{-\infty}^{0}\cap\mathcal{A}$
is non-empty (since it contains $A$), and as $T$ is $\mathbf{II}_{\infty}$, it is also non
atomic. Therefore, $\left(\xi_{-\infty}^{0}\right)^{*}\vee\mathcal{A}^{*}=\left(\xi_{-\infty}^{0}\vee\mathcal{A}\right)^{*}$
by Lemma~\ref{lem:star-partitions}, and
\[
h\left(T_{*},\left(\xi_{-\infty}^{0}\vee\mathcal{A}\right)^{*}\right)\leq H_{\mu^{*}}\left(\left(\xi_{-\infty}^{0}\vee\mathcal{A}\right)^{*}\mid T_{*}^{-1}\left(\xi_{-\infty}^{0}\vee\mathcal{A}\right)^{*}\right).\]
Since we have the other inequality, we can conclude that
$$
h\left(T_{*},\left(\xi_{-\infty}^{0}\vee\mathcal{A}\right)^{*}\right)=H_{\mu^{*}}\left(\left(\xi_{-\infty}^{0}\vee\mathcal{A}\right)^{*}\mid T_{*}^{-1}\left(\xi_{-\infty}^{0}\vee\mathcal{A}\right)^{*}\right)<\infty
$$
thus $\xi_{-\infty}^{0}\vee\mathcal{A}$ is ED.

Using this preliminary result, by considering an increasing sequence
$(\xi_{k})$ of finite local partitions with core $A$ such that $\widehat{\xi_{k}}\uparrow\mathcal{B}$,
we build an increasing sequence of ED partitions $\left(\left(\xi_{k}\right)_{-\infty}^{0}\vee\mathcal{A}\right)$.

By definition, each $\left(\left(\xi_{k}\right)_{-\infty}^{0}\vee\mathcal{A}\right)^{*}$ is a perfect $\sigma$-algebra for the corresponding factor.
In particular, $\bigcap_n T_*^{-n} \left( \left(\xi_k\right)_{-\infty}^0 \vee\mathcal{A} \right)^*$ is a zero-entropy factor of $T_*$.
As $h\left(T_{*},\left(\left(\xi_k\right)_{-\infty}^{0}\vee\mathcal{A}\right)^{*}\right)$ is finite for all $k$, we can inductively define a sequence $\eta_k = \eta_{k-1} \vee T^{-n_k}_*\left((\xi_k)^0_{-\infty}\vee\mathcal{A}\right)^*$ where the integers $\{n_k\}$ are chosen so that
$$H_{\mu^*}(\eta_i \mid \ (\eta_{j-1})^{-1}_{-\infty}) - H_{\mu^*}(\eta_i \mid \ (\eta_j)^{-1}_{-\infty})< \frac{1}{i}2^{j-i}$$ whenever $i< j$. Proceeding as in page $69$ of~\cite{parry_topics}, we obtain that $\eta:= \bigvee_{k \ge 1}\eta_k$ is a perfect $\sigma$-algebra for $T_*$.
We have to show that $\eta$ is indeed a Poissonian $\sigma$-algebra.
Observe that $T^{-n_k}_*\left((\xi_k)^0_{-\infty}\vee\mathcal{A}\right)^* =\left(T^{-n_k}\left( (\xi_k)^0_{-\infty}\vee\mathcal{A}\right) \right)^*$ and thus
$$
\eta_k = \bigvee_{j=0}^{k}\left(T^{-n_j}\left( (\xi_j)^0_{-\infty} \vee\mathcal{A}\right) \right)^*.
$$
For any $k \ge 0$, $T^{-n_k}\left( (\xi_0)^0_{-\infty}\vee\mathcal{A}\right)$ is non-atomic and for any $j \le k$, $T^{-n_k}\left( (\xi_0)^0_{-\infty}\vee\mathcal{A}\right) \subset T^{-n_j}\left( (\xi_j)^0_{-\infty}\vee\mathcal{A}\right)$. We can apply Lemma~\ref{lem:star-partitions} to get
$$
\eta_k = \left(\bigvee_{j=0}^{k}T^{-n_j}\left( (\xi_j)^0_{-\infty}\vee\mathcal{A}\right) \right)^*.
$$
Setting $\alpha_k = \bigvee_{j=0}^{k}T^{-n_j}\left( (\xi_j)^0_{-\infty}\vee\mathcal{A}\right)$ and $\alpha = \bigvee_{k \ge 1} \alpha_k$, then $\alpha_k \uparrow \alpha$ and so
$\eta_k = \alpha_k^* \uparrow \alpha^*$. We conclude that $\eta= \alpha^*$, and so $\alpha$ is a generating ED partition.

Using the fact that $\eta$ is a generating perfect $\sigma$-algebra, and applying Lemma~\ref{lem:star_monotone_algebras_roy}, we obtain
$$\mathcal{P}(T_*)= \bigcap_{n=0}^{\infty} T_*^{-n} \eta  = \left(\bigcap_{n=0}^{\infty} T^{-n} \alpha\right)^*.$$
In case  $\bigcap_{n=0}^{\infty} T^{-n} \alpha$ is $\sigma$-finite, $\mathcal{P}(T_*)$ is indeed Poissonian. Otherwise, by ergodicity of $T$, $\bigcap_{n=0}^{\infty} T^{-n} \alpha$ contains no set of positive finite measure, and $T$ is remotely infinite.
\end{proof}
%  ***************The non-ergodic case:********************
We remark that whenever $(X,\mathcal{B},\mu,T)$ is of type $\mathbf{II}_{\infty}$, but not necessarily ergodic,
by the ergodic decomposition we can uniquely decompose $\mu= \mu_0 + \mu_1$ where $\mu_1$ and $\mu_2$ are mutually singular and both are $T$-invariant,
$(X,\mathcal{B},\mu_0,T)$ is remotely infinite and $(X,\mathcal{B},\mu_1,T)$ has a Poisson-Pinsker factor as above.
%$X_{2}$ is simply obtained as the measurable union of positive finite
%measure sets measurable with respect to ${\displaystyle \bigcap_{n=0}^{\infty}}T^{-n}\alpha$.
%The sets $X_{1}:=X\setminus X_{2}$ and $X_{2}$ are $T$-invariant
%and measurable with respect to $\alpha$.
%We conclude by considering the restrictions $\alpha_{1}$ and $\alpha_{2}$ of $\alpha$ to $X_{1}$
%and $X_{2}$ respectively. Indeed ${\displaystyle \bigcap_{n=0}^{\infty}}T_{\mid X_{1}}^{-n}\alpha_{1}$
%has only sets of zero or infinite measure, i.e. $T_{\mid X_{1}}$
%is remotely infinite.
%And as $\mathcal{P}\left(\left(T_{\mid X_{2}}\right)_{*}\right)=\mathcal{P}\left(T_{*}\right)=\left({\displaystyle \bigcap_{n=0}^{\infty}}T^{-n}\alpha\right)^{*}=\left({\displaystyle \bigcap_{n=0}^{\infty}}T_{\mid X_{2}}^{-n}\alpha_{2}\right)^{*}$,
%we can set $\mathcal{P}:={\displaystyle \bigcap_{n=0}^{\infty}}T_{\mid X_{2}}^{-n}\alpha_{2}$
%which is $\sigma$-finite with respect to $\mu_{\mid X_{2}}$ to get
%$\mathcal{P}\left(\left(T_{\mid X_{2}}\right)_{*}\right)=\mathcal{P}^{*}$.

Concluding this section, we state the following proposition, which along with Proposition~\ref{prop:PoissonPinsker} and Theorem \ref{thm:poisson_rohlin_sinai} completes the picture about Poisson-Pinkser factors:
\begin{proposition}\label{prop:zero_poisson_ED}
Let $\left(X,\mathcal{B},\mu,T\right)$ be an ergodic $\mathbf{II}_{\infty}$-system
with a zero Poisson entropy factor $\mathcal{A}$. Then there exists
a generating ED partition. In particular, $\left(X,\mathcal{B},\mu,T\right)$ possesses a Poisson-Pinsker factor $\mathcal{P}$ and $\mathcal{P}^{*}$ is the Pinsker factor of~$T_{*}$.
\end{proposition}
\begin{proof}
$\mathcal{A}$ is an ED partition which satisfies $H\left( \mathcal{A}\mid T^{-1}\mathcal{A}\right)=0$ therefore, this a direct application of Theorem~\ref{thm:poisson_rohlin_sinai}.
\end{proof}

\section{Some more results, remarks and questions}
The conclusion of Proposition~\ref{prop:zero_poisson_ED} above yields a natural question:
Does the non-triviality of the Pinsker factor of $T_{*}$ imply the existence of a Poisson-Pinsker factor for $T$?
We can only partially answer this question:
\begin{proposition}\label{pro:poisson_pinsker}
Assume $T$ is an ergodic $\mathbf{II}_{\infty}$-system which satisfies $h\left(T_{*}\right)=\hpa\left(T\right)<\infty$.
If $\mathcal{P}\left(T_{*}\right)\neq\left\{ X^{*},\emptyset\right\} $
then $T$ possesses a Poisson-Pinsker factor $\mathcal{P}$ and $\mathcal{P}^{*}=\mathcal{P}\left(T_{*}\right)$.
\end{proposition}
\begin{proof}
First we are going to show that
\begin{equation}
\label{eq:Parry_generating}
\hpa\left(T\right)=\sup\left\{ H\left(\mathcal{A}\mid T^{-1}\mathcal{A}\right), \mathcal{A}\subset T^{-1}\mathcal{A},\widehat{\mathcal{A}}=\mathcal{B}\right\}.
\end{equation}
It is based on the following observation:
Consider an increasing $\sigma$-algebra $\mathcal{A}$, such that
$\widehat{\mathcal{A}}\neq\mathcal{B}$, a set $A\in\mathcal{A}$
and a finite local partition $\xi$ such that $\widehat{\mathcal{A}\vee\xi}=\mathcal{B}$
(the existence of such a partition is ensured by Theorem~2.5 in~\cite{danilenko_rudolph07} ,
as the extension $\mathcal{B}\to\mathcal{A}$ has finite relative
Poisson entropy and therefore finite relative Krengel entropy by Proposition~\ref{prop:relative_entropy}). {F}rom earlier computations, we have
\begin{eqnarray*}
 & \frac{1}{n}H_{\mu^{*}}\left(\left(\left(\xi_{-p}^{0}\right)^{*}\vee\mathcal{A}^{*}\right)_{0}^{n}\mid\mathcal{A}^{*}\right)\\
 & \to H_{\mu^{*}}\left(\left(\xi_{-p}^{0}\right)^{*}\mid T_{*}^{-1}\left(\left(\xi_{-p}^{0}\right)^{*}\right)_{-\infty}^{0}\vee\widehat{\mathcal{A}^{*}}\right)+H_{\mu}\left(\mathcal{A}\mid T^{-1}\mathcal{A}\right)
\end{eqnarray*}
and
\begin{eqnarray*}
 & \frac{1}{n}H_{\mu^{*}}\left(\left(\left(\xi_{-p}^{0}\right)^{*}\vee\mathcal{A}^{*}\right)_{0}^{n}\mid\mathcal{A}^{*}\right)\leq\frac{1}{n}H_{\mu}\left(\left(\xi_{-p}^{0}\vee\mathcal{A}\right)_{0}^{n}\mid\mathcal{A}\right)\\
 & \to H_{\mu}\left(\xi_{-\infty}^{0}\vee\mathcal{A}\mid T^{-1}\left(\xi_{-\infty}^{0}\vee\mathcal{A}\right)\right)
\end{eqnarray*}
as $n$ tends to infinity.
Therefore $H_{\mu}\left(\mathcal{A}\mid T^{-1}\mathcal{A}\right)\leq H_{\mu}\left(\xi_{-\infty}^{0}\vee\mathcal{A}\mid T^{-1}\left(\xi_{-\infty}^{0}\vee\mathcal{A}\right)\right)$
which means that for any increasing partition, we can find a generating
partition with a greater entropy. This proves~\eqref{eq:Parry_generating}.

Let $\mathcal{A}$ be an increasing and generating $\sigma$-algebra. By Lemma~\ref{lem:star_monotone_algebras_roy} we have $\left(\widehat{\mathcal{A}}\right)^{*}=\widehat{\mathcal{A}^{*}}=\mathcal{B}^{*}$,
and as $\mathcal{P}\left(T_{*}\right)$ is included in the remote
past of every generating increasing $\sigma$-algebra, we have $\mathcal{P}\left(T_{*}\right)\subset\cap_{n>0}T_{*}^{-n}\mathcal{A}^{*}$.
Since $\mathcal{P}\left(T_{*}\right)\neq\left\{ X^{*},\emptyset\right\} $,
$\mathcal{T}:=\cap_{\mathcal{A}\subset T^{-1}\mathcal{A},\widehat{\mathcal{A}}=\mathcal{B}}\cap_{n>0}T^{-n}\mathcal{A}$
possesses at least one set of non zero finite measure and therefore is a $\sigma$-finite factor.
% Now remark that since $\mathcal{P}\left(T_{*}\right)\neq\left\{ X^{*},\emptyset\right\} $,
% for any increasing $\mathcal{A}$, $\mathcal{P}\left(T_{*}\right)\subset\cap_{n>0}T_{*}^{-n}\mathcal{A}^{*}$,
% therefore $\mathcal{T}:=\cap_{\mathcal{A}\subset T^{-1}\mathcal{A},\widehat{\mathcal{A}}=\mathcal{B}}\cap_{n>0}T^{-n}\mathcal{A}$
% is not trivial and therefore a $\sigma$-finite factor.
Now,~\eqref{eq:Parry_generating} reads
$$\hpa\left(T\right)=\sup\left\{ H\left(\mathcal{A}\mid T^{-1}\mathcal{A}\vee\mathcal{T}\right),\mathcal{A}\subset T^{-1}\mathcal{A},\widehat{\mathcal{A}}=\mathcal{B}\right\} .$$
Thanks to Proposition~\ref{prop:relative_entropy}, $\sup\left\{ H\left(\mathcal{A}\mid T^{-1}\mathcal{A}\vee\mathcal{T}\right), \mathcal{A}\subset T^{-1}\mathcal{A},\widehat{\mathcal{A}}=\mathcal{B}\right\} $
is the relative Poisson entropy of $T$ with respect to $\mathcal{T}$.
But since $h\left(T_{*}\right)=\hpa\left(T\right)<\infty$,
we deduce that $h\left(T_{*},\mathcal{T}^{*}\right)=0$. We can
now apply Proposition~\ref{prop:zero_poisson_ED}.
\end{proof}

%The following corollary is an immediate application of this result:
\begin{corollary}
Assume $(X,\mathcal{B},\mu,T)$ is of type $\mathbf{II}_{\infty}$ and $h\left(T_{*}\right)=\hpa\left(T\right)<\infty$.
Assume that $f \in L^2(X,\mu)$ is a function such that
$\sigma(\{f\circ T^{n} \}_{n \in \mathbb{Z}})=\mathcal{B}$ and $f$ has singular spectral measure.
Then  $h\left(T_{*}\right)=0$.
\end{corollary}

\begin{proof}
Since $f$ has singular measure (under $T$), so has the centered
stochastic integral $f^*$ (under $T_{*}$) (i.e the
image of $f$ under the natural isomorphism between $L^{2}\left(\mu\right)$
and the first chaos of $L^{2}\left(\mu^{*}\right)$). Therefore $f^*$
is measurable with respect to $\mathcal{P}\left(T_{*}\right)$, so we
deduce $\mathcal{P}\left(T_{*}\right)\neq\left\{ X^{*},\emptyset\right\} $.
Applying Proposition~\ref{pro:poisson_pinsker}, we get that $\mathcal{P}\left(T_{*}\right)$
is a Poissonian factor. But the smallest Poissonian factor generated
by $f^*$ is the whole $\sigma$-algebra $\mathcal{B}^{*}$
(as $\mathcal{B}$ is the factor generated by $f$) and this ends
the proof.
\end{proof}

An immediate consequence of Proposition~\ref{prop:zero_poisson_ED} is the following:
\begin{corollary}\label{cor:zero_entropy_join}
If two $\mathbf{II}_{\infty}$-transformations $(X,\mathcal{B}_i,\mu_i,T_i)$  have zero Poisson entropy for $i=1,2$, then so does any joining of them.
\end{corollary}
\begin{proof}
Let $(X,\mathcal{B},\nu,T)$ be a joinning of $\mu_1$ and $\mu_2$. By Proposition~\ref{prop:zero_poisson_ED}, $(X^*,\mathcal{B}^*,\nu^*,T_*)$ has a Poissonian Pinsker factor, which contains (the pullbacks of) the $\sigma$-algebras $\mathcal{B}_i^*$ for $i=1,2$. The smallest Poissonian $\sigma$-algebra which contains these is $(\mathcal{B}_1 \vee \mathcal{B}_2)^*=\mathcal{B}^*$. Thus $(X,\mathcal{B},\nu,T)$ is its own  Poisson-Pinsker factor.
\end{proof}

Except for those cases where the Krengel and Parry entropies coincide with Poisson entropy, we do not know whether the statement corresponding to Corollary~\ref{cor:zero_entropy_join} holds with Krengel or Parry entropy.

We now state a ``strong disjointness'' result:

\begin{proposition}
If $\left(X,\mathcal{B},\mu,T\right)$ has a zero Poisson entropy factor and $\left(Y,\mathcal{C},\nu,S\right)$
has not, then they are strongly disjoint.
\end{proposition}
\begin{proof}
By Proposition~\ref{prop:zero_poisson_ED}, the Pinsker factor of $T_*$ is Poissonian.
In the proof of Proposition~\ref{prop:PoissonPinsker}, we proved that in this situation any factor of $T$ has a $\sigma$-finite intersection with its Poisson-Pinsker factor, thus has a zero Poisson entropy factor.
\end{proof}

Having already used results about relative entropy from~\cite{danilenko_rudolph07} in previous sections,
we formulate another couple of results about Poisson suspensions related to this paper of Danilenko and Rudolph:

\begin{proposition}
Let $\left(X,\mathcal{B},\mu,T\right)$ be an ergodic $\mathbf{II}_{\infty}$-system with a  Poisson-Pinkser factor $\mathcal{P}$. Then
$T$ is relatively $CPE$ (complete positive entropy) and therefore relatively weakly mixing over $\mathcal{P}$.
\end{proposition}
\begin{proof}
It is a consequence of the existence of a relative Pinsker factor with respect to a factor $\mathcal{A}$ (see Definition~1.5 in~\cite{danilenko_rudolph07}), which is the maximum factor such that any extension with respect to $\mathcal{A}$ has zero Krengel entropy. Assume $T$ admits $\mathcal{P}\left(T\right)$ as Poisson-Pinsker factor.
{F}rom Proposition~\ref{prop:relative_entropy}, relative Poisson and Krengel entropies coincide and thus give the same relative Pinsker factor. But this means that the relative Pinsker factor over $\mathcal{P}\left(T\right)$ is $\mathcal{P}\left(T\right)$ itself. This proves $T$ is relatively $CPE$ over $\mathcal{P}\left(T\right)$.
As the maximum distal extension has zero Krengel (and then Poisson) relative entropy, it is also contained in $\mathcal{P}\left(T\right)$. Thanks to the infinite Furstenberg decomposition (Proposition~4.2 in~\cite{danilenko_rudolph07}), $T$ is relatively weakly mixing over $\mathcal{P}\left(T\right)$.
\end{proof}

In \cite{danilenko_rudolph07} it is proved that a probability-preserving transformation $S$ is distal if and only if $T\times S$ is a zero entropy extension of $T$, whenever $T$ is a conservative measure-preserving transformation.
Translating this result into the Poisson framework yields the following criterion for distality:
\begin{proposition}
A probability-preserving transformation $S$ is distal if and only if, for any conservative measure-preserving transformation $T$, if the Poisson suspension $T_{*}$
has zero entropy, then so does the Poisson suspension $\left(T\times S\right)_{*}$.
\end{proposition}

We now apply our previous results to a question of Aaronson and Park from~\cite{aaro_park}, about the existence of a Krengel-Pinsker factor for
quasi-finite transformations.
First, we note that the assumptions of Theorem~\ref{thm:poisson_rohlin_sinai} hold in particular for quasi-finite systems:

\begin{corollary}\label{cor:quasi_finite_pinsker}
Let $T$ be an ergodic quasi-finite system $(X,\mathcal{B},\mu,T)$. Either it is remotely infinite or there exists a Poisson-Pinsker factor, which is also a Parry and Krengel-Pinsker factor.
\end{corollary}
\begin{proof}
Let $A \in \mathcal{B}$ be a quasi-finite sweep-out set and let $\xi:=\left\{ A,X\setminus A\right\} $ be the local quasi-finite partition induced by $A$.
Observe that the left-hand side of~\eqref{eq:return_time} is bounded by the entropy of the return times partition on $A$, thus
$H\left( \xi^0_{-\infty}\mid T^{-1}\xi^0_{-\infty} \right)$ is finite.
It follows from~\eqref{eq:EDpartition} and Lemma~\ref{lem:perfect_partitions}
%%the proof of Theorem~\ref{thm:poisson_eq_krengel_quasi_finite}
that $\xi^0_{-\infty}$ is ED.
The assumptions of Theorem~\ref{thm:poisson_rohlin_sinai} are thus satisfied.

If there is a Poisson-Pinsker factor $\mathcal{P}$, it is also a Parry-Pinsker factor by Proposition~\ref{prop:PoissonPinsker}.
To prove it is also a Krengel-Pinsker factor, we have to prove there exists a zero Krengel entropy factor.
But if $\mathcal{A}$ is the factor generated by $\xi$, then Poisson and Krengel entropy coincide on this factor and are finite.
Moreover, as in the proof of Proposition~\ref{prop:PoissonPinsker}, $\mathcal{A} \cap \mathcal{P}$ is $\sigma$-finite, therefore we can consider the extension $\mathcal{A} \cap \mathcal{P}$ to $\mathcal{A}$ where relative Poisson and Krengel entropies coincide thanks to Proposition~\ref{prop:relative_entropy}.
Since $\mathcal{A} \cap \mathcal{P}$ has zero Poisson entropy, it is also the case for Krengel entropy and we are done.
\end{proof}

%********************** Change!
Corollary~\ref{cor:quasi_finite_pinsker} generalizes a result of
from~\cite{aaro_park} about the existence of a Krengel-Pinsker factor for a special class of quasi-finite systems called $LLB$.
We do not know if the conclusion of this corollary is true without the assumption that $T$ is quasi-finite.

Since we do not know that a factor of a quasi-finite system is itself
quasi-finite, we cannot conclude that the remotely-infinite property is inherited
by factors in the quasi-finite ergodic case. However, this is the case for  $LLB$ systems which are shown to be $LLB$ on any of their factors (see~\cite{aaro_park} again). As a consequence of Theorem~\ref{thm:poisson_rohlin_sinai}, we get:

\begin{corollary}
If $T$ is $LLB$ and remotely infinite, then any factor $S$ of $T$ is remotely infinite.
\end{corollary}

The main open question left at this point, as stated in the beginning, is the following:
Are Krengel, Parry and Poisson entropies equal for \emph{every} conservative measure-preserving transformation?

At this time, we cannot answer even the following questions:
Is there an inequality between Poisson entropy and Krengel entropy which holds in general?
Are the properties of having zero Poisson entropy and having zero Krengel entropy equivalent?

Related to this is the following question from~\cite{danilenko_rudolph07}: Does any conservative transformation have a factor with arbitrarily small Poisson/Krengel entropy?
A positive answer to Danilenko and Rudolph's question would imply a positive answer to our main question. However we do not even know if there always exists a factor with \emph{finite} Poisson or Krengel entropy.

\bibliographystyle{abbrv}
\bibliography{poisson_markov}
\end{document}